\documentclass[11pt]{article}

\usepackage{amsfonts, amssymb}
\usepackage{amsbsy}
\usepackage{amsmath}
\usepackage{amsthm}
\usepackage{latexsym}
\usepackage{verbatim}
\usepackage{url}
\usepackage{graphicx,color}
\usepackage{mathrsfs}
\usepackage[normalem]{ulem}

\def\t{\mbox{\textbf{\textsf{t}}}}

\newcommand{\NN}{\mathbb{N}}

\newcommand{\ZZ}{\mathbb{Z}}

\newcommand{\TT}{\mathbb{T}}

\def\Alph{\textrm{Alph}}

\newenvironment{debut}[1]%
{\begin{quote} \bf #1 \it}%
{\end{quote}}

{\begin{debut} {R\'esum\'e~: }}
{\end{debut}}

{\begin{debut} {Abstract:}}
{\end{debut}}

\theoremstyle{plain}
\newtheorem{theo}{Theorem}[section]

\newtheorem{lem}[theo]{Lemma}
\newtheorem{cor}[theo]{Corollary}
\newtheorem{prop}[theo]{Proposition}
\newtheorem{defi}[theo]{Definition}

\theoremstyle{definition}
\newtheorem{remq}[theo]{Remark}
\newenvironment{rem}{\begin{remq}\rm }{\end{remq}}
\newtheorem{exa}[theo]{Example}
\newtheorem*{notation}{Notation}

\theoremstyle{remark}
\newtheorem*{note}{Note}

\def\t1{\mathbf{t}^{(1)}}
\def\bt{\mathbf{t}}
\def\bm{\mathbf{m}}
\def\br{\mathbf{r}}
\def\m1{\mathbf{m}^{(1)}}

\def\bs{\mathbf{s}}
\def\cA{\mathcal{A}}
\def\bbf{\mathbf{f}}

\def\bx{\mathbf{x}}

\def\bw{\mathbf{w}}
\def\rev{\widetilde}
\def\cAstar{\mathcal{A}^*}
\def\cAplus{\mathcal{A}^+}
\def\cAw{\mathcal{A}^\omega}
\def\cAinf{\mathcal{A}^\infty}
\def\cB{\mathcal{B}}

\def\PAL{\mathrm{PAL}}
\def\empt{\varepsilon}
\def\Alphit{\mbox{{\em Alph}}}
\def\TT{\mathrm{T}}
\def\cP{\mathcal{P}}
\def\cC{\mathcal{C}}

\advance \topmargin by -\headheight
\advance \topmargin by -\headsep
\numberwithin{equation}{section}

\begin{document}
\title{PALINDROMIC RICHNESS}

\author{Amy Glen\footnotemark[1] \and Jacques Justin\footnotemark[2] \and Steve Widmer\footnotemark[3] \and Luca Q.~Zamboni\footnotemark[4]}   
\date{}

\maketitle  

\footnotetext[1]{{\bf Corresponding author:} LaCIM, Universit\'e du Qu\'ebec \`a Montr\'eal, C.P. 8888, succursale Centre-ville,
Montr\'eal, Qu\'ebec, CANADA, H3C 3P8 (\texttt{amy.glen@gmail.com}). Supported by CRM, ISM, and LaCIM.}
\footnotetext[2]{LIAFA, Universit\'e Paris  Diderot - Paris 7,  
Case 7014, 75205 Paris Cedex~13, FRANCE  
(\texttt{jacjustin@free.fr}).}
\footnotetext[3]{Department of Mathematics, University of North Texas, P.O. Box  311430, Denton, Texas, TX 76203-1430, USA ({\tt sbw0061@unt.edu}).}
\footnotetext[4]{Institut Camille Jordan, 
Universit\'e Claude Bernard Lyon 1, 
43 boulevard du 11 novembre 1918, 
69622 Villeurbanne Cedex
FRANCE ({\tt luca@unt.edu}).}

\medskip
\hrule
\begin{abstract} 
In this paper, we study combinatorial and structural properties of a new class of finite and infinite words that are `rich' in palindromes in the utmost sense.  A characteristic property of so-called {\em rich words} is that all {\em complete returns} to any palindromic factor are themselves palindromes. These words encompass the well-known {\em episturmian words}, originally introduced by the second author together with X.~Droubay and G.~Pirillo in 2001. Other examples of rich words have appeared in many different contexts. Here we present the first unified approach to the study of this intriguing family of words.

Amongst our main results, we give an explicit description of the periodic rich infinite words and show that the recurrent {\em balanced} rich infinite words coincide with the balanced {\em episturmian words}. We also consider two wider classes of infinite words, namely {\em weakly rich words} and {\em almost rich words} (both strictly contain all rich words, but neither one is contained in the other). In particular, we classify all recurrent balanced weakly rich words. As a consequence, we show that any such word on at least three letters is necessarily episturmian; hence weakly rich words obey {\em Fraenkel's conjecture}. Likewise, we prove that a certain class of almost rich words obeys Fraenkel's conjecture by showing that the recurrent balanced ones are episturmian or contain at least two distinct letters with the same frequency.  

Lastly, we study the action of morphisms on (almost) rich words with particular interest in morphisms that preserve (almost) richness. Such morphisms belong to the class of {\em $P$-morphisms} that was introduced by A.~Hof, O.~Knill, and B.~Simon in 1995. \medskip

\noindent {\bf Keywords:} palindrome; episturmian word; balanced word; rich word;  return word; morphism.

\noindent MSC (2000): 68R15. 
\end{abstract}
\hrule

\section{Introduction}

In recent years there has been growing interest in palindromes in the field of {\em combinatorics on words}, especially since the work of A.~de~Luca~\cite{aD97stur} and also X.~Droubay and G.~Pirillo \cite{xDgP99pali}, who showed that the well-known {\em Sturmian words} are characterized by their {\em palindromic complexity} \cite{jAmBjCdD03pali, pBzMeP07fact, sBsHmNcR04onth}.  
A strong motivation for the study of palindromes, and in particular infinite words containing arbitrarily long palindromes, stems from their application to the modelling of {\em quasicrystals} in theoretical physics (see for instance \cite{dDlZ03comb, aHoKbS95sing}) and also Diophantine approximation (e.g., see \cite{sF06pali2}).

In \cite{xDjJgP01epis}, the second author together with X.~Droubay and G.~Pirillo observed that any finite word $w$ of {\em length} $|w|$  contains at most $|w| + 1$ distinct palindromes (including the empty word). Even further, they proved that a word $w$ contains exactly $|w| + 1$ distinct palindromes if and only if the longest palindromic suffix of any prefix $p$ of $w$ occurs exactly once in $p$ (i.e., every prefix of $w$ has {\em Property $Ju$} \cite{xDjJgP01epis}). Such words  are `rich' in palindromes in the sense that they contain the maximum number of different palindromic factors. Accordingly, we say that a finite word $w$ is {\em rich} if it contains exactly $|w| + 1$ distinct palindromes (or equivalently, if every prefix of $w$ has {\em Property $Ju$}). Naturally, an infinite word is rich if all of its factors are rich. In independent work, P.~Ambro{\v{z}}, C.~Frougny, Z.~Mas{\'a}kov{\'a}, E.~Pelantov{\'a}~\cite{pAzMePcF06pali}  have considered the same class of words which they call {\em full words}, following earlier work of S.~Brlek, S.~Hamel, M.~Nivat, and C.~Reutenauer in \cite{sBsHmNcR04onth}.

In \cite{xDjJgP01epis}, the second author together with X.~Droubay and G. Pirillo showed that the family of {\em episturmian words} \cite{xDjJgP01epis, jJgP02epis}, which includes the well-known {\em Sturmian words}, comprises a special class of  rich infinite words. Specifically, they proved that if an infinite word $\bw$ is episturmian, then any factor $u$ of $\bw$ contains exactly $|u| + 1$ distinct palindromic factors. (See \cite{jB07stur, aGjJ07epis, mL02alge} for recent surveys on the theory of Sturmian and episturmian words.) Another special class of rich words consists of S.~Fischler's sequences with ``abundant palindromic prefixes'', which were introduced and studied in~\cite{sF06pali} in relation to Diophantine approximation (see also~\cite{sF06pali2}). Other examples of rich words have appeared in many different contexts; they include the {\em complementation-symmetric sequences}~\cite{jAmBjCdD03pali}, certain words associated with $\beta$-expansions where $\beta$ is a {\em simple Parry number}~\cite{pAzMePcF06pali}, and a class of words coding $r$-interval exchange transformations~\cite{pBzMeP07fact}.

In this paper we present the first study of rich words as a whole. Firstly, in Section~\ref{S:properties}, we prove several fundamental properties of rich words; in particular, we show that rich words are characterized by the property that all {\em complete returns} to any palindromic factor are palindromes (Theorem~\ref{p1}). We also give a more explicit description of periodic rich infinite words in Section~\ref{S:periodic}~(Theorem~\ref{T:periodic}). 

In Section~\ref{S:other} we define {\em almost rich words}: they are infinite words for which only a finite number of prefixes do not satisfy {\em Property $Ju$}. Such words can also be defined in terms of the {\em defect} of a finite word $w$, which is the difference between $|w| + 1$ and the number of distinct palindromic factors of $w$ (see the work of Brlek {\it et al.} in \cite{sBsHmNcR04onth} where periodic infinite words with bounded defect are characterized). With this concept, rich words are those with defect $0$ and almost rich words are infinite words with bounded defect. `Defective words' and related notions are studied in Section~\ref{S:other}, where we also introduce the family of {\em weakly rich words} (which includes all rich words), defined as infinite words with the property that all complete returns to letters are palindromes.

In Section~\ref{S:balance} we consider applications to the {\em balance property}: an infinite word over a finite alphabet $\cA$ is {\em balanced} if, for any two factors $u$, $v$ of the same length, the number of $x$'s in each of $u$ and $v$ differs by at most $1$ for each letter $x \in \cA$. (Sturmian words are exactly the aperiodic balanced infinite words on two letters.) First we describe the recurrent {\em balanced} rich infinite words: they are precisely the balanced episturmian words. We then go much further by classifying all recurrent balanced weakly rich words. As a corollary to our classification, we show that any such word on at least three letters is necessarily episturmian. 
Consequently, weakly rich words obey {\em Fraenkel's conjecture} \cite{aF73comp}. We also prove that a certain class of almost rich words obeys Fraenkel's conjecture by showing that the recurrent balanced ones are episturmian or contain at least two distinct letters with the same frequency.

Lastly, in Section~\ref{S:morphisms}, we study the action of morphisms on (almost) rich words  with particular interest in morphisms that preserve (almost) richness. Such morphisms belong to the class of {\em $P$-morphisms} that was introduced by A.~Hof, O.~Knill, and B.~Simon in~\cite{aHoKbS95sing}  (see also the nice survey on palindromic complexity by J.-P.~Allouche, M.~Baake, J.~Cassaigne, and D.~Damanik~\cite{jAmBjCdD03pali}).

\subsection{Notation and terminology} \label{S:prelim}

In what follows, $\cA$ denotes a finite {\em alphabet}, i.e., a finite non-empty set of symbols called {\em letters}. A finite \emph{word} over $\cA$ is a finite sequence of letters from $\cA$. The {\em empty word} $\empt$ is the empty sequence. Under the operation of concatenation, the set $\cA^*$ of all finite words over $\cA$ is a {\em free monoid} with identity element $\empt$ and set of generators $\cA$. The set of {\em non-empty} words over $\cA$ is the {\em free semigroup} $\cA^+ := \cAstar \setminus \{\empt\}$.  

A (right) \emph{infinite word} $\bx$ is a sequence indexed by $\NN^+$ with values in $\cA$, i.e., $\bx = x_1x_2x_3\cdots$ with each $x_i \in \cA$. The set of all infinite words over $\cA$ is denoted by $\cAw$, and we define $\cAinf := \cAstar \cup \cAw$.  An {\em ultimately periodic} infinite word can be written as $uv^\omega = uvvv\cdots$, for some $u$, $v \in \cAstar$, $v\ne \empt$. If $u = \empt$, then such a word is {\em periodic}. An infinite word that is not ultimately periodic is said to be {\em aperiodic}.  For easier reading, infinite words are hereafter typed in boldface to distinguish them from finite words.

Given a finite word $w = x_{1}x_{2}\cdots x_{m} \in \cAplus$ with each $x_{i} \in \cA$, the \emph{length} of $w$, denoted by $|w|$, is equal to $m$ and we denote by $|w|_a$ the number of occurrences of a letter $a$ in $w$. By convention, the empy word is the unique word of length $0$.   
We denote by $\tilde w$ the {\em reversal} of $w$, given by $\tilde w = x_m \cdots x_2x_1$. If $w = \tilde w$, then $w$ is called a {\em palindrome}. The empty word $\empt$ is assumed to be a palindrome.

A finite word $z$ is a {\em factor} of a finite or infinite word $w \in \cA^\infty$ if $w = uzv$ for some $u \in \cAstar$, $v \in \cA^\infty$. In the special case $u = \empt$ (resp.~$v = \empt$), we call $z$ a {\em prefix} (resp.~{\em suffix}) of $w$. The set of all factors of  $w$ is denoted by $F(w)$ and the {\em alphabet} of $w$ is $\Alph (w) := F(w) \cap \cA$. We say that $F(w)$ is {\em closed under reversal} if for any $u \in F(w)$, $\tilde u \in F(w)$. When $w = ps \in \cA^+$, we often use the notation $p^{-1}w$ (resp.~$ws^{-1}$) to indicate the removal of the prefix $p$ (resp.~suffix $s$) of the word~$w$.

The {\em palindromic (right-)closure} of a word $u$ is the (unique) shortest palindrome $u^{(+)}$ having $u$ as a prefix  \cite{aD97stur}. That is, $u^{(+)} = uv^{-1}\tilde{u}$, where $v$ is the longest palindromic suffix of $u$. 

Given an infinite word $\bx = x_1x_2x_3\cdots$, the {\em shift map} $\TT$ is defined by $\TT(\bx) = (x_{i+1})_{i\geq 1}$ and its $k$-th iteration is denoted by $\TT^k$.  For finite words $w \in \cA^+$, $\TT$ acts circularly, i.e., if $w = xv$ where $x \in \cA$, then $\TT(w) = vx$. The circular shifts $\TT^k(w)$ with $1 \leq k \leq |w| -1$ are called \emph{conjugates} of $w$. A finite word is \emph{primitive} if it is different from all of its conjugates (equivalently, if it is not a power of a shorter word).

A factor of an infinite word $\bx$ is \emph{recurrent} in $\bx$ if it occurs infinitely often in $\bx$, and $\bx$ itself is said to be \emph{recurrent} if all of its factors are recurrent in it. Furthermore, $\bx$ is \emph{uniformly recurrent} if for all $n$ there exists a number $K(n)$ such that any factor of length at least $K(n)$ contains all factors of length $n$ in $\bx$ (equivalently, if any factor of $\bx$ occurs infinitely many times in $\bx$ with bounded gaps~\cite{eCgH73sequ}).

Let $\cA$, $\cB$ be two finite alphabets. A {\em morphism} $\varphi$ of $\cA^*$ into $\cB^*$ is a map $\varphi: \cA^* \rightarrow \cB^*$ such that $\varphi(uv) = \varphi(u)\varphi(v)$ for any words $u$, $v$ over $\cA$.  A {\em morphism on $\cA$} is a morphism from $\cA^*$ into itself. A morphism is entirely defined by the images of letters. All morphisms considered in this paper will be {\it non-erasing}, so that the image of any non-empty word is never empty. Hence the action of a morphism $\varphi$ on $\cA^*$ naturally extends to infinite words; that is, if $\bx = x_1x_2x_3 \cdots \in \cAw$, then $\varphi(\bx) = \varphi(x_1)\varphi(x_2)\varphi(x_3)\cdots$.  An infinite word $\bx$ can therefore be a {\em fixed point} of a morphism $\varphi$, i.e., $\varphi(\bx) = \bx$.  If $\varphi$ is a (non-erasing) morphism such that 
$\varphi(a) = aw$ for some letter $a \in \cA$ and $w \in \cA^{+}$,
then $\varphi$ is said to be \emph{prolongable} on $a$. In this case,
the word $\varphi^{n}(a)$ is a proper prefix of the word
$\varphi^{n+1}(a)$ for each $n \in \NN$, and the limit of the sequence
$(\varphi^n(a))_{n\geq 0}$ is the unique infinite word: 
\[
  \bw = \underset{n \rightarrow\infty}{\lim}\varphi^n(a) = \varphi^{\omega}(a) ~(= aw\varphi(w)\varphi^2(w)\varphi^3(w)\cdots).
\]
Clearly, $\bw$ is a fixed point of $\varphi$ and we say that $\bw$ is \emph{generated} by $\varphi$.

A morphism $\varphi$ on $\cA$ is said to be {\em primitive} if there exists a positive integer $k$ such that, for all $x\in \cA$, $\varphi^k(x)$ contains all of the letters of $\cA$.  Any prolongable primitive morphism generates a uniformly recurrent infinite word \cite{mQ87subs}.

For other basic notions and concepts in combinatorics on words, see for instance the Lothaire books \cite{mL83comb, mL02alge}. 

\section{Definitions and basic properties} \label{S:properties}

In this section, we prove several fundamental properties of rich words. First we recall a number of facts already mentioned in the introduction. 

\begin{prop} {\em \cite[Prop.~2]{xDjJgP01epis}}
A word $w$ has at most $|w|+1$ distinct palindromic factors.
\end{prop}

\begin{defi}
A word $w$ is {\em rich} if it has exactly $|w|+1$ distinct palindromic factors.
\end{defi}

\begin{defi}
A factor $u$ of a word $w$ is said to be {\em unioccurrent} in $w$ if $u$ has exactly one occurrence in $w$.
 \end{defi}
 
 \begin{prop}{\em \cite[Prop.~3]{xDjJgP01epis}} \label{P:xDjJgP01epis} A word $w$ is rich if and only if all of its prefixes (resp.~suffixes) have a unioccurrent palindromic suffix (resp.~prefix). 
\end{prop}

\begin{cor}
If $w$ is rich, then:
\begin{enumerate}
\item[i)] it has exactly one {\em unioccurrent palindromic suffix} (or \emph{ups} for short);
\item[ii)] all of its factors are rich;
\item[iii)] its reversal $\tilde w$ is also rich. \qed
\end{enumerate} 
\end{cor}

\begin{note} $i)$ is {\em Property $Ju$} from \cite{xDjJgP01epis}.
\end{note}

Clearly, if $w$ has a ups, $u$ say, then $u$ is the {\em only} ups of $w$, and moreover $u$ is the longest palindromic suffix of $w$. So if $w = vu$, then $w^{(+)} = vu\tilde v$. Furthermore:

 \begin{prop} \label{P:palindromic-closure}
 Palindromic closure preserves richness.
 \end{prop}
 \begin{proof}
 Let $w$ be rich with ups $u$. The case $w=u$ is trivial. Now suppose $w=fu$ for some (non-empty)  word $f = f'x$, $x \in \cA$. Then $w^{(+)}=fu \tilde f = f' xux \tilde f'$. Clearly $xux$ is a ups of $fux$, so  continuing we see that all prefixes of $fu \tilde f$ have a ups of the form $hu \tilde h$ where $h$ is a suffix of $f$. Thus $w^{(+)}=fu \tilde f$ is rich.
\end{proof}

\begin{prop}
If $w$ and $w'$ are rich with the same set of palindromic factors, then they are {\em abelianly equivalent}, i.e., $|w|_x = |w'|_x$ for all letters $x \in \cA$.
\end{prop}
\begin{proof} Any palindromic factor of $w$ (resp.~$w'$) ending (and hence beginning) with a letter $x \in \cA$ is the ups of some prefix of $w$ (resp.~$w'$). Thus the number of $x$'s in $w$ (resp.~$w'$) is the number of palindromic factors ending with $x$. 
\end{proof}

\begin{prop}
Suppose $w$ is a rich word. Then there exist letters $x$, $z\in\Alphit(w)$ such that $wx$ and $zw$ are rich.
\end{prop}
\begin{proof}
If $w$ is rich and not a palindrome, then let $u$ be its ups with $|u|<|w|$. Then $u$ is preceded by some letter $x$ in $w$, thus $wx$ has $xux$ as its ups, and hence $wx$ is rich. If, on the contrary, $w$ is a palindrome, let $w=xv$ and let $u$ be the ups of $v$. If $u=v$ then $v$ is a palindrome so $w=xv$ gives $v=x^n$ for some $n$ whence $wx=x^{n+2}$ which is rich. If $|u|<|v|$ then let $y$ be the letter before $u$ in $v$. Then $yuy$ is ups of $vy$ (because $u$ is the ups of $v$). If $yuy$ is not the ups of $wy$ then it is a prefix of $w$, but then $u$ is a prefix of $v$, and as $|u|<|v|$, $u$ occurs twice in $v$, contradiction. Thus $yuy$ is the ups of $wy$ and this one is rich. 

In view of Proposition~\ref{P:xDjJgP01epis}, one can similarly show that $zw$ is rich for some letter $z \in \Alph(w)$. 
\end{proof}  
\begin{note} In the case when $w$ is rich and not a palindrome, the fact that $wx$ is rich for some letter $x \in \Alph(w)$ is a direct consequence of Proposition~\ref{P:palindromic-closure}.
\end{note}

Naturally:

\begin{defi}
An infinite word  is {\em rich} if all of its factors are rich.
\end{defi}

\begin{prop} There exist recurrent rich infinite words that are not uniformly recurrent. 
\end{prop}
\begin{proof}
Consider the infinite word $\bt$ generated by the morphism $(a \mapsto aba,\ b \mapsto bb)$ from \cite{jC97comp}. It suffices to show (rather easily) that $\bt$ is rich. Similarly, the {\em Cantor word} of \cite{nP02subs}, or even $abab^2abab^3abab^2abab^4abab^2abab^3abab^2abab^5\cdots$ (fixed point of the morphism: $a \mapsto abab$, $b \mapsto b$)  are recurrent rich infinite words that are not uniformly recurrent. 
\end{proof}
 
 \begin{prop} \label{P:reverse-closure}
 A rich infinite word $\bs$ is recurrent if and only if its set of factors $F(\bs)$ is closed under reversal.
 \end{prop}
 
 The proof of this proposition uses the following lemma. Note that ``richness'' is not necessary for the ``if'' part.
 
 \begin{lem} \label{L:palindromic-prefixes}
 A recurrent rich infinite word has infinitely many palindromic prefixes.
  \end{lem}
 \begin{proof} 
 Let $v_1$ be a non-empty prefix of a recurrent rich infinite word $\bs$. Being rich, $v_1$ has a unioccurrent palindromic prefix, $u_1$ say (by Proposition \ref{P:xDjJgP01epis}). Let $v_2$ be a prefix of  $\bs$ containing a second occurrence of $u_1$. It has a unioccurrent palindromic prefix, $u_2$ say. Now, $u_2$ is not a prefix of $u_1$ because $u_1$ is not unioccurrent in $v_2$, thus $|u_2|>|u_1|$. 
 \end{proof}
 
 \begin{rem} Although the well-known {\em Thue-Morse word} $\bm$, which is the fixed point of the morphism $\mu: a\mapsto ab, b \mapsto ba$ beginning with $a$, contains arbitrarily long palindromes  (see Example~\ref{ex:Thue-Morse} later), $\bm$ is not rich. For instance, the prefix $abbabaabba$ is not rich (since its longest palindromic suffix $abba$ is not unioccurrent in it).
 \end{rem}
 
 \begin{proof}[Proof of Proposition $\ref{P:reverse-closure}$]
IF: Consider some occurrence of a factor $u$ in $\bs$ and let $v$ be a prefix of $\bs$ containing $u$. As $F(\bs)$ is closed under reversal, $\tilde v \in F(\bs)$. Thus, if $v$ is long enough, there is an occurrence of  $\tilde u$ strictly on the right of this particular occurrence of $u$ in $\bs$. Similarly $u$ occurs on the right of this $\tilde u$ and thus $u$ is recurrent in $\bs$. 

ONLY IF: As $\bs$ is recurrent, it follows from Lemma \ref{L:palindromic-prefixes} that $F(\bs)$ is closed under reversal.
  \end{proof}
 
\newpage
\begin{theo}\label{p1}
For any finite or infinite word $w$, the following properties are equivalent:
 \begin{itemize}
\item[i)] $w$ is rich;
 
\item[ii)] for any factor $u$ of $w$, if $u$ contains exactly two occurrences of a palindrome $p$ as a prefix and as a suffix only, then $u$ is itself a palindrome.
\end{itemize}
 \end{theo}

 \begin{proof} 
 $i) \Rightarrow ii)$: Suppose, on the contrary, $ii)$ does not hold for rich $w$. Then $w$ contains a non-palindromic factor $u$ having exactly two occurrences of a palindrome $p$ as a prefix and as a suffix only. Moreover, these two occurrences of $p$ in $u$ cannot overlap. Otherwise $u = pv^{-1}p$ for some word $v$ such that $p = vf = gv = \tilde v\tilde g = \tilde p$; whence $v = \tilde v$ and $u = g\tilde v \tilde g = gv\tilde g$, a palindrome. So $u=pzp$ where $z$ is a non-palindromic word. We easily see that $u$ does not have a ups; thus $u$ is not rich, a contradiction.

$ii) \Rightarrow i)$: Otherwise, let $u$ be a factor of $w$ of minimal length satisfying $ii)$ and not rich. Trivially $|u|>2$, so let $u=xvy$ with $x, y \in \cA$. Then $xv$ is rich by the minimality of $u$. Since $u$ is not rich whilst $xv$ is rich, the longest palindromic suffix $p$ of $u$ occurs more than once in $u$. Hence, by $ii)$ we reach a contradiction to the maximality of $p$. 
 \end{proof}
 
 \begin{remq} \label{R:crw}
 Given a finite or infinite word $w$ and a factor $u$ of $w$, we say that a factor $r$ of $w$ is a {\em complete return} to $u$ in $w$ if $r$ contains exactly two occurrences of $u$, one as a prefix and one as a suffix ({\em cf.} `first returns' in \cite{cHlZ99desc}). With this notion, Property $ii)$ says that all complete returns to any palindromic factor are themselves palindromes.  
In particular, consecutive occurrences of a letter $x$ in a rich word are separated by palindromes.  
\end{remq}

\begin{note} In view of Theorem~\ref{p1}, an alternative proof of the richness of episturmian words  can be found in the paper \cite{vAlZiZ03pali} where the fourth author, together with V.~Anne and I.~Zorca, proved that for episturmian words,  all complete returns to palindromes are palindromes. See also \cite{jJlV00retu} for further work on `return words' in Sturmian and episturmian words.
\end{note}

 \section {Periodic rich infinite words} \label{S:periodic}

Theorem~\ref{p1} provides a characterization of rich infinite words by complete returns to palindromes. We now give a more explicit description of periodic rich infinite words.

 \begin{theo}  \label{T:periodic}
For a finite word $w$, the following properties are equivalent:
 \begin{itemize}
\item[i)] $w^\omega$ is rich;  
\item[ii)] $w^2$ is rich; 
 \item[iii)] $w$ is a product of two palindromes and all of the conjugates of $w$ (including itself) are rich.  
 \end{itemize}
 \end{theo}

\begin{exa}
$(aabbaabab)^\omega$ and $(abcba)^\omega$ are rich.
\end{exa}

The proof of Theorem~\ref{T:periodic} requires several lemmas. In what follows, $x$ and $z$ always denote letters.

 \begin{lem} \label{l1}
 If $u$ is rich and $ux$ has a palindromic suffix $r$ such that $2|r| \ge |u|$, then $ux$ is rich.
  \end{lem}
 \begin{proof}
 We can suppose $r$ has maximal length. If $r$ has another occurrence in $ux$, then, as  $2 |r|+1 \ge |ux|$, the two occurrences overlap or are separated by at most one letter. Thus they both form a palindrome which  is a suffix of $ux$ and is strictly longer than $r$, a contradiction. Therefore $r$ is the ups of $ux$,  which is rich. 
 \end{proof}

\begin{lem} \label{l2}
 If $w=pq$, $p,q$ palindromes, then $w$ has a conjugate $w'=p'q'$, $p',q'$ palindromes with $||p'|-|q'|| \le 2$.
 \end{lem}
 \begin{proof}
 Easy.
 \end{proof}
 
 \begin{lem} \label{l3}
 If $w=pq$, $p,q$ palindromes, is rich and $2|q| \ge |w|-4$ (resp.~$2|p| \ge |w|-4$), then $pqp$ (resp.~$qpq$) is rich.
  \end{lem}
 \begin{proof}
 Suppose $2|q| \ge |w|-4$ (the other case is obtained by reversal as $\tilde w$ is rich).
 If $pqp$ is not rich, let $vz$, $v \in \cAstar$, be the shorter prefix of $p$ such that $pqvz$ is not rich. Further, let $r$ be the longest palindromic suffix of $pqvz$. Then, as $z \tilde v q vz$ is a suffix of $pqvz$,  we have $|r| \ge |q|+2|v| +2$; whence $2 |r| \ge 2|q|+4|v| +4 \ge |w| + 4|v|\ge |w|$. Then by Lemma \ref{l1} $pqvz$ is rich, contradiction. 
 \end{proof}
 
\begin{lem} \label{l4}
 If $w=pq$, $p,q$ palindromes, and $pqp,qpq$ are rich, then $w^2$ is rich.
  \end{lem}
 \begin{proof}
 If $w^2=pqpq$ is not rich, let $vz$ be the shorter prefix of $q$ such that $pqpvz$ is not rich. As $qpq$ is rich, its prefix $qpvz$ has a ups, $r$ say. As $z \tilde v pvz$ is a palindromic suffix of $qpvz$, $r$ must begin in the prefix $q$ of $qpvz$. As $r$ is not the ups of $pqpvz$, consider its leftmost occurrence in $pqpvz$. If the two occurrences overlap or are separated by at most one letter, both they form a palindromic suffix of $pqpvz$. As this one is not the ups of $pqpvz$, it has another occurrence; whence $r$ has another occurrence on the left of the leftmost one, contradiction.
 
 Thus, the two considered occurrences of $r$ do not overlap. This implies that the leftmost occurrence of $r$ lies in the prefix $pq$ of $pqpvz$, but then by reversal $r$ also occurs in $qp$, hence $r$ is not the ups of $qpvz$, contradiction.    
 \end{proof}

\begin{proof}[Proof of Theorem $\ref{T:periodic}$] $i) \Rightarrow ii)$: Trivial.\par

$ii) \Rightarrow iii)$:  It suffices to show that $w$ is a product of two palindromes.
Let $r$ be the ups of $w^2$. Then, clearly $|r| > |w|$, thus $r = qw$ and $w=pq$ for some $p,q$. Therefore $r=qpq=\tilde q \tilde p \tilde q$, whence $p$ and $q$ are palindromes.\par

$ii) \Rightarrow i)$:  We show first that $w^3$ is rich. By $iii)$, $w^2$ has a ups $qw$ and $w = pq$, $p,q$ palindromes. For any $u,v$ such that $uv=p$, consider $f=w^2 u$. Observe that $f$ has a palindromic suffix $\tilde u qpq u$ which is its ups, otherwise $qpq$ would not be the ups of $w^2$. Thus all such $f$ are rich, in particular $w^2p$ is rich.
Now, if $ez$ is a prefix of $q$, we show by induction on $|e|$ that $w^2pez$ is rich. Let $r$ be the longest palindromic suffix of  $w^2pez$. As this one has suffix $z \tilde e pqp ez$, we have $|r|\ge 2|p|+|q|+2|e| +2$; whence $2|r|  \ge |w^2pez|$. Thus, by Lemma \ref{l1}, $w^2pez$ is rich, and hence $w^3$ is rich.

Now denote by $\bs _n$ the prefix of length $n$ of $\bs = w^\omega$. We show by induction on $n$ that $w^3 \bs_n$ is rich. Let $r$ be the ups of $w^2 \bs_n$. Then $r$ is also a suffix of $w^3 \bs_n$. Clearly these two occurrences of $r$ overlap, thus both give a palindromic suffix of $w^3 \bs_n$. If this one were not the ups of $w^3 \bs_n$, there would be another occurrence of $r$ in $w^2 \bs_n$, contradiction. Thus $w^3 \bs_n$ has a ups and, as $w^3 \bs_{n-1}$ is rich, it is rich too.\par

$iii) \Rightarrow i)$: By Lemma \ref{l2}, $w$ has a (rich) conjugate $w'=p'q'$ with $p', q'$ palindromes and $||p'|-|q'|| \le 2$, whence by Lemma \ref{l3} $p'q'p'$ and $q'p'q'$ are rich. Thus by Lemma \ref{l4} $(w')^2$ is rich. So, using part ``$ii) \Rightarrow i)$'', $(w')^\omega$ is rich, and so too is $w^\omega$.
\end{proof}

\begin{rem}
For $iii)$ the hypothesis that all of the conjugates of $w$ are rich is not sufficient: $abc$ is so, but $(abc)^\omega$ is not rich. The hypothesis that $w$ is rich and a product of two palindromes is not sufficient: $w=ba^2bab^2aba^2b$ is a rich palindrome, but $\TT(w) = a^2bab^2aba^2b^2$ is not rich.  
\end{rem}

\section{Some related words} \label{S:other}

\subsection{Defects \& oddities}

The {\em defect}  \cite{sBsHmNcR04onth}  of a finite word $w$ is defined by
\[
  D(w) = |w| + 1 - |\PAL(w)|,
\]  
where $\PAL(w)$ denotes the set of distinct palindromic factors of $w$ (including~$\empt$). This definition naturally extends to infinite words $\bw \in \cAw$ by setting $D(\bw)$ equal to the maximum defect of the factors of $\bw$. In fact, this definition may be refined by observing that if $u$ is a factor of a word $v$, then $D(u) \leq D(v)$ (see \cite{sBsHmNcR04onth}); thus
\[
  D(\bw) = \max\{D(u) ~|~ \mbox{$u$ is a prefix of $\bw$}\}.
\]  

With this notion, finite or infinite rich words are exactly those with defect equal to $0$ (called {\em full words} in \cite{pAzMePcF06pali, sBsHmNcR04onth}). Accordingly, we  say that an infinite word with bounded defect is {\em almost rich}. Such infinite words contain only a finite number of prefixes that do not have a ups.

\begin{notation}
Let $\bt_n$ denote the prefix of length $n$ of a given finite or infinite word $\bt$.
\end{notation}

\begin{prop} \label{P:defect-positions}
If $\bt_n$ has a ups, then $D(\bt_n)= D(\bt_{n-1})$, otherwise  $D(\bt_n)= D(\bt_{n-1})+1$.
\end{prop}
\begin{proof} 
If $\bt_n$ has a ups, then $\bt_n$ contains one more palindromic factor than $\bt_{n-1}$, whence $D(\bt_n) = D(\bt_{n-1})$. On the other hand, if $\bt_n$ has no ups, then $\bt_n$ has the same number of palindromic factors as $\bt_{n-1}$, thus $D(\bt_{n}) = D(\bt_{n-1}) + 1$.
\end{proof}

In other words, if $\bt$ has defect $k$, then there are exactly  $k$ ``defective" positions; hence it is appropriate to say that such a word $\bt$ has $k$ defects.

\begin{remq} A noteworthy fact is that for a given word $w$ with $k$ defects, the extension $wx$, with $x \in \Alph(w)$, may not have the same number of defects (in particular, the palindromic closure of $w$ may have greater defect). For example, $w=caca^2bca$ has 2 defects,  but $wx$ has $3$ defects for $x = a$, $b$, or $c$. 
\end{remq}
Periodic almost rich words have the following simple characterization.

\begin{theo}\label{T:per-alm}
A periodic infinite word $w^\omega$ is almost rich if and only if $w$ is a product of two palindromes.
\end{theo}
\begin{proof} 
The ``if'' part follows immediately from \cite[Theorem 6]{sBsHmNcR04onth}:  {\em if $p$, $q$ are palindromes and $pq$ is a primitive word, then the defect of $(pq)^\omega$ is bounded by the defect of its prefix of length $|pq| + \lfloor \frac{||p| - |q||}{3} \rfloor$}. Conversely, if $w^\omega$ is almost rich, then, for large enough $n$, $w^n$ has a ups. Thus, as in the proof of $ii) \Rightarrow iii)$ of Theorem~\ref{T:periodic}, we get $w= q p $ for some palindromes $p$, $q$. 
\end{proof}

\begin{prop} \label{P:bounded-reversal}
If an almost rich word $\bt$ is recurrent, then $F(\bt)$ is closed under reversal.
\end{prop}
\begin{proof}  Let $u$ be any prefix of $\bt$ with the same defect number $k$ as $\bt$. By recurrence, we can consider another occurrence of $u$ such that $\bt= su \cdots$ for some non-empty word $s$. Then, any suffix $s'u$ of $su$ has a ups since every prefix $v$ of $\bt$ with $|v| > |u|$ has a ups (otherwise the defect of $\bt$ would be greater than $k$, by Proposition \ref{P:defect-positions}). In particular, $su$ has a ups, $p$ say. Now, $p$ is not a suffix of $u$ because $u$ is not unioccurrent in $su$, so $|p| > |u|$ and we have $p = vu = \tilde{u}\tilde{v}$ for some non-empty word $v$. Thus $u$ and $\tilde{u}$ are both factors of $\bt$, and hence $F(\bt)$ is closed under reversal. 
\end{proof}

\begin{defi} \label{D:oddity}
The pair $\{w, \tilde w \}$ is an {\em oddity} of a finite or infinite word $\bt$ if either $w$ or $\tilde w$ (or both!) is a non-palindromic complete return to some non-empty palindromic factor of $\bt$ (called the {\em incriminated palindrome}).
 \end{defi}

\begin{note} An oddity of a finite or infinite word $\bt$ takes the form $pup$ where $p$ is the incriminated palindrome and $u$ is a non-palindromic word. Indeed, non-palindromic complete returns to any palindrome $p$ are necessarily longer than $2|p|+1$ (see the proof of Theorem~\ref{p1}).
\end{note}

Let $O(\bt)$ denote the number of oddities of $\bt$.

\begin{prop} \label{P:oddity} 
$O(\bt) \leq D(\bt)$.
\end{prop}
\begin{proof}
Let $w=pup$ be an oddity of $\bt$ with $p$ the incriminated palindrome. Let $n$ be the minimal integer such that $w$ or $\tilde w$ occurs in $\bt_n$ (thus as a suffix). If $\bt_n$ has a ups, $q$ say, then $|q| > |p|$, trivial. If  $|q| < |w|$, then  the prefix $p$ of $q$ occurs in the interior of the complete return $w=pup$, impossible. If $|q| > |w|$ then $\tilde w$ occurs as a prefix of $q$, in contradiction with minimality of $n$. Thus $\bt_n$ does not have a ups, i.e., $n$ is a defective position. Thus each oddity gives a defect. For achieving the proof we have to show that $n$ cannot be the same for two different oddities. Suppose $\bt_n$  has a second (suffix) oddity $ w'=qvq$. If $|p|=|q|$, clearly $w'=w$. Otherwise let $|q|<|p|$ for instance, then $q$ occurs twice in $p$, thus $w'$ is a suffix of $p$, hence $\tilde w'$ is a prefix of $p$, contradicting the minimality of $n$.        
\end{proof}

\begin{exa} 
We may have $O(\bt) < D(\bt)$; for instance the periodic word $(abcabcacbacb)^\omega$ has 3 oddities ($abca$, $bcab$, $cabc$) ending at positions $4$, $5$, $6$, but 4 defects at positions $4$, $5$, $6$, $7$.
The periodic infinite word $(abc)^\omega$ has a defect at each position $n \ge 4$ but only three oddities ($abca$, $bcab$, $cabc$). So infinitely many defects do not necessarily give rise to infinitely many oddities.
\end{exa}

\begin{prop} \label{P:infinite-oddities}
A uniformly recurrent infinite word has infinitely many oddities if and only if it has infinitely many palindromic factors and infinitely many defects.
\end{prop}

\begin{proof} ONLY IF: Suppose $\bs$ is a uniformly recurrent infinite word with infinitely many oddities. Clearly, $\bs$ has infinitely many defects as $D(\bs) \geq O(\bs)$ by Proposition~\ref{P:oddity}. Moreover, $\bs$ must have infinitely many palindromic factors. Otherwise, if $\bs$ contains only a finite number of palindromes, then each of its palindromic factors has only a finite number of different return words (and hence only a finite number of non-palindromic complete returns) as $\bs$ is uniformly recurrent.  Hence $\bs$ has only finitely many oddities, a contradiction. 

IF: Suppose, by way of contradiction, $\bs$ is a (uniformly recurrent) infinite word with infinitely many palindromic factors and infinitely many defects, but only a  finite number of oddities. Then there are only finitely many palindromic factors of $\bs$ that are incriminated by the oddities and the longest of these palindromes has length $L$, say. Since $\bs$ contains infinitely many palindromes, it has infinitely many `non-defective' positions. Thus there exists an arbitrarily large  $n$ such that $n$ is a defective position and such that the prefix $\bs_{n+1} =\bs_n x$, with $x \in \cA$, has a ups $q$, with $|q| > L + 2$. Thus, with $q = x q' x$, we see that $q'$ is a palindromic suffix of $\bs_{n}$. Let $r$ be the longest palindromic suffix of $\bs_{n}$. As $n$ is a defective position, $r$ has another occurrence in $\bs_{n}$ and $rur$ is a non-palindromic complete return to $r$, which is a contradiction since $|r| \geq |q'| > L$. 
\end{proof}

\begin{remq} Uniform recurrence is necessary for the ``only if'' part of the above proposition (but useless for the ``if'' part). For example, with $v_1 = abcd$ and $v_{n} = v_{n-1}(abc)^n d$ for $n\geq 2$, the (non-uniformly) recurrent  infinite word $v_1 v_2 v_3 \cdots$ has infinitely many oddities, but only five palindromic factors: $\empt$, $a$, $b$, $c$, $d$.
\end{remq}

\begin{exa} \label{ex:Thue-Morse} The Thue-Morse word $\bm$ has infinitely many oddities. Indeed, since $\bm$ is generated by the morphism $\mu^2: a\mapsto abba, b \mapsto baab$, $\bm$ clearly contains infinitely many palindromes. Moreover, one can prove by induction that $\bm$ has infinitely many  defects occurring in runs of length $2^{2n+1}$ starting at positions  $2^{2n+3} +1$ and $2^{2n+4} + 2^{2n+3} + 1$ for $n\geq 0$.
\end{exa}

\subsection{Weakly rich words}

We say that an infinite word $\bw$ over $\cA$ is {\em weakly rich} (or simply a {\em WR-word} for short) if for every $a\in\cA$,  all complete returns to $a$ in $\bw$ are palindromes. This class of words contains all rich words but is in fact a much larger class. Clearly every binary word is weakly rich but not necessarily rich. The periodic infinite word $(aacbccbcacbc)^\omega$ is readily verified to be weakly rich but not rich (since the complete return to $aa$ is not a palindrome). Note, however, that the family of weakly rich words neither contains nor is contained in the family of almost rich words. Indeed, the WR-word  $(aacbccbcacbc)^\omega$ has infinite defect (since it does not take the form $(pq)^\omega$ with $p$, $q$ palindromes, and hence contains only finitely many distinct palindromic factors -- see Theorem~\ref{T:per-alm}). There also exist almost rich words that are not weakly rich; for instance, the almost rich word $(aabacabaac)^\omega$  (which has only $2$ defects at positions $10$ and $11$) is not weakly rich as $cabaac$ is a non-palindromic complete return to $c$.

Our motivation for introducing WR-words will become evident in the next section. 


\section{Applications to balance} \label{S:balance}

A finite or infinite word is said to be {\em balanced} if, for any two of its factors $u$, $v$ with $|u| = |v|$, we have $||u|_x-|v|_x|\leq 1$ for any letter $x$, i.e., the number of $x$'s in each of $u$ and $v$ differs by at most~$1$. {\em Sturmian words} are precisely the aperiodic balanced infinite words on a $2$-letter alphabet.

{\em Fraenkel's conjecture} \cite{aF73comp} is  a well-known problem related to balance that arose in a number-theoretic context and has remained unsolved for over thirty years. Fraenkel conjectured that, for a fixed $k \geq 3$, there is only one covering of $\ZZ$ by $k$ {\em Beatty sequences} of the form $(\lfloor \alpha n + \beta \rfloor)_{n\geq1}$, where $\alpha$, $\beta$ are real numbers. A combinatorial interpretation of this conjecture may be stated as follows (taken from  \cite{gPlV07acha}). Over a $k$-letter alphabet with $k\geq 3$, there is only one recurrent balanced infinite word, up to letter permutation  and shifts, that has mutually distinct letter frequencies. This supposedly unique infinite word is called {\em Fraenkel's sequence} and is given by $(F_k)^\omega$ where the {\em Fraenkel words} $(F_i)_{i\geq 1}$ are defined recursively by $F_1=1$ and  $F_{i}=F_{i-1}iF_{i-1}$ for  all $i\geq 2.$ 
For further details, see \cite{gPlV07acha} and references therein.

In \cite{gPlV07acha}, Paquin and Vuillon characterized balanced episturmian words by classifying these words into three families. Amongst these classes, only one has mutually distinct letter frequencies and, up to letter permutation and shifts, corresponds to Fraenkel's sequence. That is:

\begin{prop} {\em \cite{gPlV07acha}} \label{P:gPlV07acha} Suppose $\bt$ is a balanced episturmian word with $\Alphit(\bt) = \{1,2,\ldots, k\}$, $k \geq 3$. If $\bt$ has mutually distinct letter frequencies, then up to letter permutation, $\bt$ is a shift of~$(F_k)^\omega$.
\end{prop}

In this section, we first show  that recurrent balanced rich infinite words are necessarily (balanced) episturmian words. Then, using a special map, we classify all recurrent balanced weakly rich words. As a corollary, we show that any such word (on at least three letters) is necessarily a (balanced) episturmian word. Thus, although WR-words constitute a larger class of words than episturmian words, the subset of those which are balanced coincides with those given by Paquin and Vuillon in \cite{gPlV07acha}. Consequently, WR-words obey {\em Fraenkel's conjecture}. Using techniques similar to those in the rich case, we also prove that a certain class of almost rich words (with only a few oddities) obeys Fraenkel's conjecture by showing that the recurrent balanced ones are episturmian or contain at least two distinct letters with the same frequency.

Before proceeding, let us recall some useful well-known facts about balance (see for instance the survey~\cite{lV03bala} and references therein). 
\begin{itemize}
\item In a balanced word, the gaps between successive occurrences of any letter $x$ belong to a pair $\left\{ k,k+1 \right\}$ for some integer $k\geq 0$. 
\item Any recurrent balanced infinite word with alphabet $\cA$ and $|\cA| > 2$ is periodic.
\end{itemize}

\subsection {Balanced rich words} \label{S:balance-rich}

The main result of this section is the following theorem.

\begin{theo} \label{T:balance}
Recurrent balanced rich infinite words are precisely the balanced episturmian words.
\end{theo}

First we prove some lemmas. For a given letter $a$, we denote by $\psi_a$ the morphism defined by $\psi_a: a\mapsto a, x \mapsto ax$ for all letters $x \ne a$. A noteworthy property of $\psi_a$ is that $\psi_a(w)a$ is a palindrome if and only if $w$ is a palindrome. It is well-known that an infinite word is {\em epistandard} if and only if it is generated by an infinite composition of the morphisms $\psi _x$. Moreover, an infinite word $\bt$ is {\em episturmian} if and only if $F(\bt) = F(\bs)$ for some epistandard word $\bs$ (see \cite{xDjJgP01epis, jJgP02epis}).

\begin{lem}\label{l5}
Suppose $\bs=\psi_a(\bt)$ for some letter $a$ and infinite word $\bt$.
\begin{itemize}
\item[i)] If $\bs$ is rich, then $\bt$ is rich.
\item[ii)] If $\bs$ is balanced, then $\bt$ is balanced.
\end{itemize}
\end{lem}
\begin{proof}
$i)$: If false, let $wx \in F(\bt)$ be minimal such that $wx$ is not rich. If $x \neq a$ then let $r$ be the ups of $\psi_a(wx)$ which is rich by hypothesis. Then $ar=\psi_a(h)$ where $h$ is a palindromic suffix (but not a ups) of $wx$. Thus $h$ has another occurrence in $wx$, which implies $r$ has another occurrence in  $\psi_a(wx)$, contradiction.

If $x=a$ we consider the ups of $\psi_a(wa)a=\psi_a(w)aa$ and by a similar argument we reach a contradiction.  \medskip

$ii)$: If $\bt$ is not balanced, then it contains two factors $u$, $v$ of the same minimal length such that $||u|_x - |v|_x|=2$ for some $x$. Let $U=\psi_a(u)$, $V=\psi_a(v)$, then $|U|=2|u|-|u|_a$ and $|V|=2|v|-|v|_a$. By adding to and/or deleting from $U,V$ some $a$ we get $U',V'$ factors of $\bs$ of the same length. If $x \neq a$ then $||U'|_x - |V'|_x|=2$. If $x=a$ then, as $|U|_a=|u|= |v| = |V|_a$, we get  $||U'|_a - |V'|_a|=2$. In both cases, $\bs$ is not balanced, contradiction.
\end{proof}

\begin{remq} If $\bs = \psi_a(\bt)$ or $\bs = a^{-1}\psi_a(\bt)$ for some letter $a$ and infinite word $\bt$, then the letter $a$ is \emph{separating for $\bs$} and its factors; that is, any factor of $\bs$ of length 2 contains the letter $a$.
\end{remq}

\begin{lem} \label{L:separating}
Suppose $\bt$ is a recurrent infinite word with separating letter $a$ and first letter $x \neq a$. Then $\bt$ and   $a \bt$ have the same set of factors.
\end{lem}

\begin{proof} Clearly $F(\bt) \subseteq F(a\bt)$. To show that $F(a\bt) \subseteq F(\bt)$, let $u$ be any factor of $a\bt$. If $u = a$ or $u$ is not a prefix of $a\bt$,  then clearly $u \in F(\bt)$. Otherwise, if $u \ne a$ is a prefix of $a\bt$, then $u$ takes the form $axu'$ where $xu'$ is a prefix of $\bt$. As $\bt$ is recurrent, $xu'$ occurs again in $\bt$, and hence $u = axu'$ must be a factor of $\bt$ because the letter  $a$ is separating for~$\bt$. 
\end{proof}

This almost trivial lemma allows us to ignore the cases where the separating letter is not the first letter.

\begin{proof}[Proof of Theorem~$\ref{T:balance}$]
Let $\bs$ be a rich, recurrent and balanced infinite word. If  $\bs$ has a separating letter, $a$ say, then $\bs=\psi_a(\bt)$ or $\bs = a^{-1}\psi_a(\bt)$ for some recurrent infinite word $\bt$, which is also rich and balanced  by Lemmas~\ref{l5} and \ref{L:separating}. If we can continue infinitely in this way then $\bs$ is episturmian by the work in \cite{jJgP02epis}. Otherwise we arrive at some recurrent infinite word, rich and balanced, without a separating letter; call it $\bt$. In particular, no $xx$ occurs in $\bt$ (because $x$ would be separating). We call such an infinite word without factor $xx$, $x \in \cA$, a {\it skeleton}.  Consider any factor of form $xpx$ of $\bt$ with $p$ $x$-free. By Theorem~\ref{p1}, $p$ is a palindrome, and as no square of a letter occurs in it, $p$ has odd length. If $xp_1x$ and $xp_2x$ 
are two such factors of $\bt$, then in view of balance, $|p_1|=|p_2|$. Thus $x$ occurs in $\bt$  with period $\pi_x=|p_i|+1$. Take for $x \in \cA$ the letter with minimal $\pi_x$ (i.e., the letter with the greatest frequency in $\bt$). If $y$ is any other letter, as $ \pi_y \ge \pi_x$, only one $y$ may occur in a  $xp_ix$. By symmetry, this $y$ lies at the centre of $p_i$. Thus all $p_i$ are reduced to their centre, i.e., $\pi_x=2$ and $x$ is separating, contradiction. So this case is impossible and $\bs$ is episturmian.          
\end{proof}

As an immediate consequence of Proposition~\ref{P:gPlV07acha} and Theorem~\ref{T:balance}, we have:

\begin{cor} \label{C:rich-balanced}
Recurrent balanced rich infinite words with mutually distinct letter frequencies are Sturmian words or have the form given by Fraenkel's conjecture. \qed
\end{cor}

\subsection{Balanced weakly rich words}

We now establish a much stronger result, namely that WR-words obey Fraenkel's conjecture, by proving that balanced WR-words on at least three letters are necessarily (balanced) episturmian words. First we classify all such words. In order to state our classification, we need the following notation.

Let $\bx = x_1x_2x_3\cdots \in \cA^\omega$ with each $x_i \in \cA$, and let $a$ be a new symbol not in $\cA$.  We define $ \sigma_a: \cA^\omega \rightarrow (\cA\cup \{a\})^\omega$ by
\[
\sigma_a (\bx)= ax_1a^{\epsilon_1}x_2a^{\epsilon_2}x_3a^{\epsilon_3}\ldots
\]
where $\epsilon _i \in \{1,2\},$ with $\epsilon _i=2$ if and only if $x_i = x_{i+1}$.

\begin{theo}\label{T:weakly-rich} Suppose $\bw$ is a recurrent balanced WR-word with $\Alphit(\bw) =\{1,2,\ldots ,k\}$, $k \geq 3$. Then, up to letter permutation, $\bw$ is either:
\begin{itemize}
\item[$1)$] a shift
of the periodic word
\[\psi_1^n\circ \psi_2\circ \cdots \circ \psi_{k-1}(k^\omega) \quad \mbox{for some $n\geq 1;$}\] 

\item[$2)$] or a shift of the periodic word
\[\sigma_1\circ \sigma_2\circ \cdots \circ \sigma_j\circ \psi_{j+1}^2\circ \psi_{j+2}\circ \cdots \circ \psi_{k-1}(k^\omega) \quad \mbox{for some $1\leq j\leq k-2$.}
\] 
\end{itemize}
\end{theo}

The proof of Theorem~\ref{T:weakly-rich} requires several lemmas. In what follows, we assume that $\Alph(\bw) = \cA$ with $|\cA|\geq 3$. For each $a\in \cA$, we set $g_a=\sup |u|$ where the supremum is taken over all factors $u$ of $\bw$ not containing the letter $a$. 

First we recall a useful lemma from \cite{Tij1}.

\begin{lem} {\em \cite[Lemma 6]{Tij1}} \label{L:tijdeman}
Suppose $\bw \in \cAw$ is balanced, and let $a \in\cA$ be such that the frequency of $a$ in $\bw$ is at least $1/3$. Then the word  $\bw'\in (\cA\setminus\{a\})^\omega$ obtained from $\bw$ by deleting all occurrences of the letter $a$ in $\bw$ is also balanced.
\end{lem}

\begin{lem}\label{bal} Suppose $\bw \in \cA^\omega$ is a recurrent balanced WR-word, and let $a\in \cA$ be such that $g_a \leq g_x$ for all $x\in \cA.$ Then the word $\bw'\in (\cA\setminus\{a\})^\omega$ obtained from $\bw$ by deleting all occurrences of the letter $a$ in $\bw$ is also a recurrent  balanced WR-word.
\end{lem}

\begin{proof} Clearly $\bw'$ is a recurrent WR-word; in fact, for each letter $x\neq a,$ every complete return to $x$ in $\bw '$ is a complete return to $x$ in $\bw$ with all occurrences of $a$ deleted. It remains to show that $\bw '$ is balanced. Since $\bw$ is balanced, it follows that if $aUa$ is a complete return to $a$ in $\bw,$   then each $x\in \cA$ occurs at most once in $U.$ Otherwise, if some letter $x$ occurred more than once in $U,$ we would have $g_x<g_a.$ Moreover, since $\bw$ is a WR-word, $U$ must be a palindrome. Thus $|U|\leq 1$, and hence $g_a=1.$ It follows that the frequency of $a$ in $\bw$ is at least equal to  $1/3,$ and hence from Lemma~\ref{L:tijdeman}, we deduce that the word $\bw '$ obtained from $\bw$ by deleting all occurrences of $a$ is  balanced. 
\end{proof}

\begin{lem}\label{sigma}  Let $\bw $ and $\bw '$ be as in Lemma~$\ref{bal}$. Suppose $\bw '$ contains the factor $bb$ for some $b\in \cA\setminus\{a\}.$ Then $\bw $ is a shift of $\sigma_a(\bw ').$ In particular, the complete returns to $b$ in $\bw$ are of the form $baab$ or $baxab$ for some $x\in \cA\setminus\{a,b\}.$
\end{lem}

\begin{proof} Assume $\bw '$ contains the word $bb.$ Since $\bw '$ is balanced (see Lemma~\ref{bal}), every factor of $\bw'$ of length $2$ must contain at least one occurrence of $b,$ and hence  $\bw'$ contains the factor $bxb$ for every $x\in \cA\setminus\{a,b\}.$ 
Since $g_a=1$ (see Lemma~\ref{bal}), it follows that $\bw$ contains factors of the form $ba^lb$ and $ba^kxa^kb$ for every $x\in \cA\setminus\{a,b\}$ with both $l,k\geq 1.$ It is readily verified from the balance property that if  $ba^{k_1}xa^{k_1}b$ and $ba^{k_2}ya^{k_2}b$ are both factors of $\bw$ with $x,y\in \cA\setminus\{a,b\},$  then $k_1=k_2.$ 
Again by the balance property it follows that
 $|2k+1-l|\leq 1$  and $k\in \{l, l+1, l-1\}.$
If $k=l,$ then $|l+1|\leq 1$ from which it follows that $l=0, $ a contradiction.
If $k=l+1,$ then $|l+3|\leq 1,$ again a contradiction. If $k=l-1,$ then $|l-1|\leq 1$ from which it follows that $l=2$ and $k=1,$ for otherwise either $l$ or $k$ would equal $0.$ Thus $\bw$ is obtained
from $\bw '$ by inserting one $a$ between any pair of consecutive distinct letters in $\bw '$ and two $a$'s between consecutive $b$'s. In other words, $\bw$ is a shift of $\sigma _a(\bw ')$, 
as required.
\end{proof}

\begin{lem}\label{cube}   Suppose $\bw $ and $\bw '$ are as in Lemma~$\ref{bal}$ and let $b\in \cA\setminus\{a\}.$  Then $bbb$ is not a factor of $\bw '.$ 
\end{lem}
\begin{proof}Suppose $\bw '$ contains $bbb.$ Then by Lemma~\ref{sigma}, $\bw$ contains
the factors $baabaab$ and $bacab$ for some $c\in \cA\setminus\{a,b\}.$  But since $\bw$ is balanced, $\bw$ cannot contain both $aabaa$ and $bacab.$
\end{proof}

\begin{proof}[Proof of Theorem~$\ref{T:weakly-rich}$] We prove Theorem~\ref{T:weakly-rich} by induction on the number of letters $k.$
Suppose $\bw$ is a  recurrent balanced WR-word on the alphabet $\cA_3=\{1,2,3\}.$ Without loss of generality we can assume $g_1\leq g_2\leq g_3.$  Let $\bw '\in \{2,3\}^\omega$ be the word
obtained from $\bw$ by deleting all occurrences of $1$ in $\bw .$ 
First suppose $22$ does not occur in $\bw '.$ In this case $\bw '$ is a shift of the periodic word
$(23)^\omega  =\psi _2(3^\omega).$ So the only complete return to $2$ in $\bw '$ is $232.$
It follows that there exists an $n\geq 1$ such that the only complete return to $2$ in $\bw'$ is $21^n31^n2.$
Hence $\bw =\psi_1^n\circ \psi_2 (3^\omega).$ 
Next suppose $22$ occurs in $\bw '.$ It follows from the above lemmas that the complete returns to $3$ in $\bw '$ are of the form $323$ and $3223.$ But if both factors occurred in $\bw',$ then by Lemma~\ref{sigma}, $\bw$ would contain both $31213$ and $31211213,$ which contradicts the fact that
$\bw$ is balanced. Thus $\bw '$ is a shift of the periodic word $(223)^\omega =\psi_2^2(3^\omega),$ and hence by Lemma~\ref{bal}, $\bw$ is a shift of $\sigma _1\circ \psi _2^2 (3^\omega ).$
Thus, Theorem~\ref{T:weakly-rich} holds for $k=3.$

Next take $k>3$ and suppose that $\bw$ is a  recurrent balanced  WR-word on the alphabet $\cA_k=\{1,2,\ldots ,k\}.$ By induction the hypothesis we assume Theorem~\ref{T:weakly-rich} holds for any recurrent balanced WR-word on an alphabet of size $k-1.$ Without loss of generality we can assume that $g_1\leq g_2\leq \cdots \leq g_k.$ Let $\bw '$ be the word on the alphabet $\{2,3,\ldots ,k\}$ obtained from $\bw$ by deleting all occurrences of $1$ in $\bw .$ It follows from Lemma~\ref{bal} that $\bw '$ is a  recurrent balanced  WR-word, and hence by the induction hypothesis, $\bw '$ is either a shift of
$\psi_2^n\circ \psi_3\circ \cdots \circ \psi_{k-1}(k^\omega)$ for some $n\geq 1,$ or else a shift of
$\sigma_2\circ \sigma_3\circ \cdots \circ \sigma_j\circ \psi_{j+1}^2\circ \psi_{j+2}\circ \cdots \circ \psi_{k-1}(k^\omega)$ for some $2\leq j\leq k-2.$ 

First suppose that $22$ does not occur in $\bw '.$
In this case $\bw '$ must be a shift of $\psi_2\circ \psi_3\circ \cdots \circ \psi_{k-1}(k^\omega).$
Thus the complete returns to $2$ in $\bw '$ are all of the form $2x2$ for some $x\in \{3,4,\ldots ,k\}.$
Hence there exists an $n\geq 1$ such that each complete return to $2$ in $\bw$ is of the form $21^nx1^n2$
where $x\in \{3,4,\ldots ,k\}.$ Thus in this case $\bw$ is a shift  of $\psi_1^n\circ \psi_2\circ \cdots \circ \psi_{k-1}(k^\omega).$ 

Next suppose $\bw'$ contains the factor $22.$ Then by Lemma~\ref{cube}, 
$\bw '$ is either a shift of $\psi_2^2\circ \psi_3\circ \cdots \circ \psi_{k-1}(k^\omega),$ or a shift
of $\sigma_2\circ \sigma_3\circ \cdots \circ \sigma_j\circ \psi_{j+1}^2\circ \psi_{j+2}\circ \cdots \circ \psi_{k-1}(k^\omega)$ for some $2\leq j\leq k-2.$ It follows from Lemma~\ref{sigma} that
$\bw$ is either a shift of $\sigma _1 \circ \psi_2^2\circ \psi_3\circ \cdots \circ \psi_{k-1}(k^\omega),$
or else a shift of $\sigma_1\circ \sigma_2 \circ \cdots \circ \sigma_j\circ \psi_{j+1}^2\circ \psi_{j+2}\circ \cdots \circ \psi_{k-1}(k^\omega)$ for some $2\leq j\leq k-2.$ Thus $\bw$ is a shift of
$\sigma_1 \circ \sigma_2 \circ \cdots \circ \sigma_j\circ \psi_{j+1}^2\circ \psi_{j+2}\circ \cdots \circ \psi_{k-1}(k^\omega)$ for some $1\leq j\leq k-2$, as required. This concludes our proof
of Theorem~\ref{T:weakly-rich}.
\end{proof}

\begin{cor} Suppose $\bw$ be a recurrent  balanced WR-word with $\Alphit(\bw) = \{1,2,\ldots ,k\}$, $k\geq 3.$ Then $\bw$ is a (balanced) periodic episturmian word.
\end{cor}
\begin{proof} Recall that any infinite word generated by an infinite composition of the morphisms $\psi _i$ is episturmian. Thus $\psi_1^n\circ \psi_2\circ \cdots \circ \psi_{k-1}(k^\omega)$ is a periodic episturmian word. It remains to show that the words described in case 2) of Theorem~\ref{T:weakly-rich}  are periodic episturmian words. To do this, we use the Fraenkel words $(F_i)_{i\geq1}$ (defined previously). It is readily verified that if $\mathbf{x}=x_1x_2x_3\ldots$ is an infinite word not containing the symbols $\{1,2,\ldots ,n\},$ then
\begin{eqnarray}\label{eq:(2.1)} \sigma _1\circ \sigma _2\circ \cdots \circ \sigma _{j}\circ \psi _{j+1}^2(\mathbf{x})=F_{j+1}^2x_1F_{j+1}^2x_2F_{j+1}^2x_3\ldots \end{eqnarray}
and
\begin{eqnarray*}\psi _1\circ\psi _2\circ \cdots \circ \psi _{j+1}\circ \psi _1(\mathbf{x})=F_{j+1}^2x_1F_{j+1}^2x_2F_{j+1}^2x_3\ldots  \end{eqnarray*}
It follows that
\[\sigma_1\circ \sigma_2\circ \cdots \circ \sigma_j\circ \psi_{j+1}^2\circ \psi_{j+2}\circ \cdots \circ \psi_{k-1}(k^\omega )=\psi_1\circ \psi_2\circ \cdots \circ \psi_j\circ \psi_{j+1}\circ \psi _1 \circ \psi_{j+2}\circ \cdots \circ \psi_{k-1}(k^\omega ).\] Since the right hand side above is an infinite periodic episturmian word, it follows that the periodic infinite words listed in Theorem~\ref{T:weakly-rich}
are episturmian words.
\end{proof}

Hence, by Proposition~\ref{P:gPlV07acha}, WR-words obey Fraenkel's conjecture; in fact, we can show this rather easily without the use of Proposition~\ref{P:gPlV07acha}.

\begin{cor} Suppose $\bw$ is a recurrent balanced WR-word with $\Alphit(\bw) = \{1,2,\ldots ,k\}$, $k \geq 3$. If $\bw$ has mutually distinct letter frequencies,  then up to letter permutation, $\bw$ is a shift of $(F_k)^\omega.$ 
\end{cor}
\begin{proof} By Theorem~\ref{T:weakly-rich}, $\bw$ is isomorphic to a shift of one of the two types of periodic words listed in the statement of the theorem. We note that except for the extreme case of $j=k-2$ in case~$2)$ of Theorem~\ref{T:weakly-rich}, the frequency of the symbols $k$ and $k-1$ are equal.
Thus, under the added hypothesis that distinct letters occurring in $\bw$ have distinct frequencies, we deduce that $\bw$ is isomorphic to a shift of
\[\sigma_1\circ \sigma_2\circ \cdots \circ \sigma_{k-2}\circ \psi_{k-1}^2 (k^\omega).\]
By \eqref{eq:(2.1)}, we have
\[\sigma_1\circ \sigma_2\circ \cdots \circ \sigma_{k-2}\circ \psi_{k-1}^2 (k^\omega)=
(F_{k-1}^2k)^\omega = F_{k-1}(F_k^\omega) \]
which is clearly a shift of the Fraenkel sequence $(F_k)^\omega .$
\end{proof}

\subsection{Balanced almost rich words}

We now extend our study to words having only a {\em few} oddities. In the spirit of Lemma \ref{l5}, we first prove the following result (see also Theorem~\ref{T:special-defect} to follow).

\begin{prop} \label{P:defective} If $\bs = \psi_a(\bt)$, then $D(\bs) \geq D(\bt)$; in particular, if $\bs$ is  almost rich then $\bt$ is almost rich.
\end{prop}

\begin{exa}  The periodic infinite word $\bt = (abcbac)^\omega$ has $1$ defect and $\bs = \psi_a(\bt) = (a^2bacaba^2c)^\omega$ has $2$ defects. More generally, for any $k \geq 1$, $\bt = (a^kba^{k-1}ca^{k-1}ba^kc)^\omega$ has $k$ defects (see \cite{sBsHmNcR04onth}), so applying $\psi_a$ to $\bt$ gives a periodic infinite word with $k+1$ defects.
\end{exa}

\begin{proof}[Proof of Proposition $\ref{P:defective}$] If $\bs$ is rich, then $\bt$ is rich (by Lemma~\ref{l5}), and hence $D(\bs) = D(\bt) = 0$. So now suppose that $\bt$ has at least one defect. Consider any   prefix $\bt_m$ of $\bt$ corresponding to a defect, i.e., $\bt_m$ does not have a ups. Let $\bt_m=\bt_{m-1}x$ where $x$ is a letter. We show that if $x\neq a$ (resp.~$x=a$), then  $\psi_a(\bt_m)$ (resp.~$\psi_a(\bt_m)a$) has no ups and thus gives a defect in $\bs$. Let $q$ be the longest palindromic suffix of $\bt_m$ which is not unioccurrent in $\bt_m$ since $\bt_m$ has no ups. \medskip

\noindent{\em Case $x \neq a$}: $p= a^{-1} \psi_a(q)$ is the longest palindromic suffix of $\psi_a(\bt_m)$, which is not unioccurrent in it, otherwise $q$ has another occurrence in $\bt_m$, a contradiction. \medskip 

\noindent{\em Case $x=a$}: Similar to the above case, but with $\psi_a(\bt_m)a$ and its longest palindromic suffix given by $p= \psi_a(q)a$.   
\end{proof}

\begin{note} Proposition~\ref{P:defective} can be extended without difficulty to oddities; that is, if $\bs = \psi_a(\bt)$, then $O(\bs) \geq O(\bt)$.
\end{note}

The main result of this section is the following:

\begin{theo} \label{T:almost-rich-balance} Suppose $\bs$ is a recurrent balanced infinite word with alphabet $\cA$, $|\cA| >2$, and less than $|\cA|$ oddities. Then $\bs$ is either episturmian or two of its letters have the same frequency. 
\end{theo}
\begin{proof} The proof relies on  two lemmas, which are stated and proved below. As in the proof of Theorem~\ref{T:balance}, we decompose $\bs$ as much as possible using morphisms $\psi_x$, $x \in \cA$. If we can continue infinitely, then $\bs$ is episturmian. Otherwise we halt at some skeleton $\bt$ without a separating letter and with alphabet $\cB$. If $|\cB| <3$,  $\bt$ is Sturmian and hence has a separating letter, a contradiction. If $|\cB| >2$ then by Lemma \ref{L:alph} $O(\bt) < |\cB|$. Moreover, by Lemma \ref{L:skel}, as $\bt$ has no separating letter, it takes the following form (up to a shift): $\bt = (x(ab)^n ax(ba)^n b)^\omega$, for some $n \ge 1$.

If the decomposition from $\bs$ to $\bt$ uses neither $\psi_a$ nor $\psi_b$, then $a$ and $b$ have same frequencies in $\bs$ (as in $\bt$), as claimed. Otherwise we have for instance 
$\bs = \mu_1 \psi _a \mu_2 (\bt)$ with $\psi_a$, $\psi_b$ not occurring in $\mu_2$. Then considering factors $xab$ and $bab$ of $\bt$ we have $\mu_2(xab)= fxgagb$ and $\mu_2(bab)= gbgagb$ for some  $\left\{a,b \right\}$-free words $f$, $g$; whence  $\psi_a(xgagb) =ax h a hab$ and  $\psi_a(bgagb) =ab h a hab$ where $h=\psi_a(g)$, showing the unbalance $axhaha, bhahab$. Thus $\bs$ in not balanced, a contradiction.
\end{proof}

\begin{lem}
 \label{L:alph} Let $\bm = \psi_c (\br)$ and suppose $\bm$ is balanced with alphabet $\cA$, $|\cA| >2$,  and less than $|\cA|$ oddities. Then, if $\br$ has alphabet $\cB= \cA \setminus \left\{c \right\}$, $\br$ is a (balanced) skeleton with less than $|\cB|$ oddities.
\end{lem}
\begin{proof}
If $w =x_1x_2 \cdots x_n$ is an oddity in $\br$ then  $x_1cx_2 \cdots c x_n$ and $c x_1 \cdots c x_n c$ are oddities in $\bm$, thus $ 2 O(\br)=O(\bm)<|\cA|$, which implies $O(\br)<|\cB|$ if $|\cA|\ge 3$. It is also clear that if $\br$ is not a skeleton, then it contains some $aa$, whence $aca \in F( \bm)$. As $|\cA| >2$,  there is another letter $b$ in $\Alph(\bm)$, and hence $cbc \in F( \bm)$. Thus $\bm$  is not balanced, a contradiction.
\end{proof}

\begin{lem} \label{L:skel} Let $\bt$ be a recurrent balanced skeleton with alphabet  $\cB$, $|\cB| >2$, and less than  $|\cB|$ oddities and without any separating letter. Then, up to a shift, $\bt$ takes the form $(x(ab)^n ax(ba)^n b)^\omega$ for some $n \ge 1$.  
\end{lem}
\begin{proof}
As $O(\bt)< |\cB|$ there is some letter, $x$ say, such that all of the complete returns to $x$ are palindromes (of the same odd length); call them $xvx$, $xv'x$, $xv''x$, $\ldots$ and write $v=uz \tilde u$, $v'=u'z' \tilde u'$ and so on. We have $|u|>0$ (otherwise $x$ is separating in $\bt$). Consider a factor $xvxv'x$. Suppose firstly that $u=u'$. If $u$ is not a palindrome, let $u= ea \cdots b \tilde e$, $a \neq b$. Then we have factors $b \tilde ezeb$ and $a \tilde e x e a$, contradicting the balance property. Thus $u$ is an (odd) palindrome, say $u = w y \tilde w$. By the same argument, $w$ is a palindrome and so on. Thus $u$ has the form $w_n$ for some $n$, with $w_{i+1}= w_i y_i w_i$ and $w_1$, $y_i \in \cB$.

If all $u$, $u'$, $\ldots$ are equal, then the letter $w_1$ is separating in $\bt$, a contradiction.    
Thus we have for instance $u \ne u'$. Then $u= \cdots ae$, $u'= \cdots be$ with $e \in F(\bt)$ and $a$,$ b \in \cB$. The factors $aez \tilde e a$ and $bez' \tilde e b$ give $z=b$, $z'=a$. Clearly, $a$ and $b$ do not occur in $e$, otherwise we have for instance $a f b \tilde f a$ and  $a f a \tilde f a$ being factors of $\bt$ with $f$ $a$-free. So, by the gap property for $a$, $|f| \ge 2  |f|$; whence $f= \empt$. But then $aa$ is a factor of $\bt$, a contradiction.  Now observe that $u$, $u'$ have the following property.

Let $u(i)$, $u'(i)$ be the $i$-th letter of $u$, $u'$, respectively. If $u(i) \not\in \{a,b\}$, then $u'(i) = u(i)$. On the other hand, if $u(i) \in \left\{ a,b \right\}$, then $u'(i) \in \left\{a,b \right\}$. The proof is easy using $u b \tilde u$ and $u' a \tilde u'$. Thus we can write 
\[ u=f_0c_0f_1c_1 \cdots f_n c_n e, \ u'=f_0c'_0f_1c'_1 \cdots f_n c'_n e,\ f_i \in \cB \setminus \left\{ a,b \right\}, \ c_i,c'_i \in  \left\{ a,b \right\}. \]
We  easily see that $f_0= \tilde e$ and $f_i= \tilde f_{n+1-i}$, using $\tilde u x u'$ and $ub \tilde u$. Now consider $c_0 \tilde e x e c'_0$ in $\tilde u x u'$. If $c_0 = c'_0 =a$ for instance, then as $b$ does not occur in it, $b$ has a gap greater than $2 |e| +2$, a contradiction. Thus $c_0 =a$, $c'_0 =b$ for instance. But now for the same reason we have the factor $ba \tilde e x e b a$, i.e., $f_1=f_n= \empt$. It follows that $a$ and $b$ have the same gaps: $2 |e| +1$, $2 |e| +2$; moreover $a$ and $b$ alternate in $u$, $u'$. Thus 
\begin{equation} u= \tilde e ab f_2 a \cdots b a e, \ u'= \tilde e ba f_2 b \cdots a b e \label{E:odd} \end{equation}
or
\begin{equation} u= \tilde e b a f_2 a \cdots b a e, \ u'= \tilde e a b f_2 a \cdots a b e \label{E:even} . \end{equation}

Consider the first case (the second one is similar). The factor $baeb \tilde e ab$ in $ub \tilde u$ shows  that $|a e| $ is a gap for $b$, and hence $|e|+1 \ge 2|e|+1$. Thus $e= \empt$ and we have 
\begin{equation} u=  ab f_2 a \cdots b a , \ u'=  ba f_2 b \cdots a b  .\end{equation}
Now $1+ |f_i|$ is a gap for $a$ with $1+ |f_i| \le 2$; thus $f_i \in \cB$  or $f_i= \empt$ .
Moreover, considering $f_i a f_{i+1}$ for instance, we have $|f_i| + 1+|f_{i+1}| \le 2$; thus if $f_i$ is a letter, then $f_{i-1}$ and $f_{i+1}$ are empty.

Now consider $xvxv'xv''x$. If $u'=u''$ then factors $bxb$ and $aba$ contradict the balance property. Thus $u'\ne u''$ and easily  $u''=u$, $z''=z=b$; whence, up to a shift, $\bt= (xub \tilde u x u' a \tilde u')^ \omega$.

Any letter $y=f_i \in \cB \setminus \left\{ x,a,b \right \}$ gives rise to two oddities, namely $af_iba$ and $bf_i ab$, and the left-most occurrence of $y$ in $u$ gives an oddity: $ya \cdots x \cdots by$ for instance. Also $abxa$ and $baxb$ are oddities. Therefore  $O(\bt) \ge 3(|\cB| -3) +2= 3|\cB| -7$ and, as $O(\bt)< |\cB|$, this gives  $|\cB| \le 3$. Thus all $f_i$ are empty and, for some $k\geq 0$, $u=(a b)^k a$ and $v= (a b)^{2k+1} a$. Similarly, the form of equation (\ref{E:even}) gives  $v= (a b)^{2k} a$, $k \ge 1$. Hence, up to a shift, $\bt= (x(ab)^n a x(ba)^n b)^ \omega$ for some $n>0$ and letter $x$.                                                    
\end{proof}

Thus we get another class of infinite words, wider than ``rich'', that obey Fraenkel's conjecture.

\section{Action of morphisms} \label{S:morphisms}

In this section, we study the action of morphisms on (almost) rich words, with particular interest in morphisms that ``preserve'' (almost) richness. We say that a morphism $\varphi$ on $\cA$ {\em preserves} (resp.~{\em strictly preserves}) a property $P$ of (finite or infinite) words if  $w \in \cAinf$ has property $P \Rightarrow \varphi(w)$ has property $P$ (resp.~$w \in \cAinf$ has property  $P \Leftrightarrow \varphi(w)$ has property $P$).

\medskip
\begin{note} For ``richness'', finite or infinite words give the same definition for ``preserves'' (but not for ``strictly preserves''). For ``almost richness'' the definition has meaning only for infinite words.
\end{note}

\subsection{Various results} \label{SS:various}

Part $i)$ of Lemma \ref{l5}  works in the opposite sense; thus we have:

\begin{prop} \label{P:ext-l5}
Let $\bs=\psi_a(\bt)$. Then $\bs$ is rich if and only if $\bt$ is rich. 
\end{prop}
\begin{proof}
It suffices to show the ``if" part. If $\bs$ is not rich, then let $wx$ be the shortest prefix of $\bt$ such that $\psi_a(wx)$ is not rich. We show first that $\psi_a(w)a$ is rich. Let $p$ be the ups of $w$. Then  $\psi_a(w)a$ ends with the palindrome $\psi_a(p)a$. If this one has another occurrence in $\psi_a(w)a$ then $\psi_a(w)a = g \psi_a(p) a ha$, $h \in \cAstar$ whence $w=g'ph'$, $g= \psi_a(g')$, $ah= \psi_a(h')$, thus $h' \neq \empt$ and $p$ is not unioccurrent in $w$, a contradiction.

Now if $x=a$ the proof is over, otherwise it remains to show that $\psi_a(wx)=\psi_a(w)ax$ has a ups. Suppose $q$ is the ups of $wx$ (which exists since $wx$ is rich). Then $q$ begins and ends with $x \neq a$, and hence $a^{-1}\psi_a(q)$ is a palindromic suffix of $\psi_a(wx)$. As previously, we easily see that it is unioccurrent in $\psi_a(wx)$.
\end{proof}

\begin{rem} The ``if and only if'' part of Proposition~\ref{P:ext-l5} does not extend to finite words. For instance, with $v=abca$, $\psi_a(v)= aabaca$ is rich while $v$ is not rich.  
\end{rem}

\begin{cor} \label{C:epi-strict}
{\em Episturmian morphisms} strictly preserve richness of infinite words. 
\end{cor}
\begin{proof}
Proposition~\ref{P:ext-l5} and Lemma~\ref{L:separating} show that any {\em elementary epistandard morphism} $\psi_a$, as well as its {\em conjugate} $\bar{\psi}_a: a \mapsto a, x \mapsto xa$,  strictly preserve richness of infinite words. Consequently, {\em episturmian morphisms} \cite{xDjJgP01epis,jJgP02epis} strictly preserve richness of infinite words as the monoid of all such morphisms is generated by all the $\psi_a$, $\bar{\psi}_a$, and permutations of the alphabet. 
\end{proof}

\begin{prop}\label{P:insert} For a fixed letter $a \in \cA$, the `insertion' morphism $\varphi_a$, defined by $\varphi_a: x \mapsto xa$ for all $x \in \cA$, preserves richness.
\end{prop}
\begin{proof}
Let $p$ be the ups of a rich word $u$. If $p \neq u$ then $a \varphi_a(p)$ is clearly a ups of $\varphi_a(u)$, but we also have to show that $ \varphi_a(u)a^{-1}$ has a ups: this one is $ \varphi_a(p)a^{-1}$. Now if $p = u$ then $ \varphi_a(u)a^{-1}$ is its own ups. Also let  $u=yt$ with $y \in \cA$ and let $q$ be the ups of $y$. If $t$ is a palindrome, then $u =a^n$ for some $n$, a trivial case. Otherwise, let $q$ be the ups of $t$. Then $r=a \varphi_a(q)$ is the ups of $a \varphi_a(t)$ and it cannot be a prefix of $ \varphi_a(u)$ because otherwise, as $q$ is a prefix of $u$, we get $q= a^n$ for some $n$; whence easily we have a contradiction. 
\end{proof}

The next proposition deals with a transformation which is not a morphism in general. Let $w$ be a finite or infinite word. For any letter $a \in \Alph(w)$, if $a^k x$ is a prefix of $w$ (or $xa^k$ is a suffix) or $ya^k x$ occurs in $w$ with $x, y \neq a$, we say that $k$ is an {\em exponent of $a$ in $w$}. Let $ k_1<k_2 <\cdots$ be the sequence of the exponents of $a$ in $w$ and let $ h_1<h_2 <\cdots$ be another sequence of positive integers of the same length with $h_i \leq k_i$ for all $i$. Let $\pi_a(w)$ be the word obtained by replacing every exponent $k_i$ by $h_i$ in  $w$. Then:
 
\begin{prop} \label{P:pi} $\pi_a(w)$ is rich if and only if $w$ is rich. \end{prop}
\begin{proof} 
Suppose $w$ is rich. Then by Theorem~\ref{p1} the complete returns to any palindromic factor of $w$ are also palindromes. The same is true for $\pi_a(w)$ since $\pi_a$ strictly preserves palindromes (i.e., a finite word $u$ is a palindrome if and only if $\pi_a(u)$ is a palindrome). Hence $\pi_a(w)$ is rich (again by Theorem~\ref{p1}). The converse is proved similarly. 
\end{proof}

\begin{prop} \label{P:rich-fixpoint} If $\varphi$ preserves richness and is prolongable on $a \in \cA$, then $\varphi^\omega(a)$ is a rich infinite word.
\end{prop}
\begin{proof} 

  This is a trivial consequence of the fact that, for all $n\geq 1$, $\varphi^n(a)$ is a rich word,  since $\varphi(a)$ is a rich word and $\varphi$ preserves richness.
\end{proof}

\begin{note} The converse does not hold. For example, the morphism $\delta: a \mapsto aba$, $b \mapsto bcb$, $c \mapsto cbc$ generates rich infinite words, beginning with $a$, $b$, and $c$ as easily seen; however, $\delta$ does not preserve richness (e.g., $\delta(acb) = abacbcbcb$ has a defect at the second occurrence of the letter $b$). 
\end{note}

Clearly, a morphism $\varphi$ on $\cA$ preserves palindromes if and only if $\varphi(x)$ is a palindrome for all $x \in \cA$.

\begin{prop} \label{P:strict-injective} Suppose $\varphi$ is a morphism on $\cA$, with $|\cA|>1$. If $\varphi$ strictly preserves palindromes, then $\varphi$ is injective.
\end{prop}
\begin{proof} Suppose $\varphi$ strictly preserves palindromes and assume $\varphi(u)=\varphi(v)$ for some non-empty words $u$, $v \in \cA^*$. Then, with $p=u\tilde u$ and $q=v\tilde v$, $\varphi(p)=\varphi(q)$ is a palindrome. Indeed, both $\varphi(p)$ and $\varphi(q)$ are palindromes since $\varphi$ preserves palindromes, and moreover
\[
\varphi(p) = \varphi(u)\varphi(\tilde u) = \varphi(v)\rev{\varphi(u)} = \varphi(v)\rev{\varphi(v)} = \varphi(v)\varphi(\tilde v)= \varphi(q).\]
Whence $\varphi(pq) = \varphi(p)^2$ is a palindrome and $pq$ too (since $\varphi$ {\em strictly} preserves palindromes). Therefore $pq=qp$, and hence $p$ and $q$ are powers of a common word (e.g., see Lothaire~\cite{mL83comb}), i.e., $p= w^m$ and $q= w^n$. Therefore, since $\varphi(p) = \varphi(q)$, we must have $m = n$; whence $u=v$. Thus  $\varphi$ is injective.
\end{proof}

\begin{exa} \label{ex:preserve-pali} The non-injective morphism $\varphi: a \mapsto aba$, $b \mapsto bcb$, $c \mapsto aba$ preserves palindromes, but not strictly as $\varphi(abc) = ababcbaba$ is a palindrome whereas the preimage $abc$ is not.
\end{exa}

The {\em letter-doubling morphism} $\varphi_d$ defined by $\varphi_d: x \mapsto xx$ for all $x \in \cA$ strictly preserves palindromes; it also preserves almost richness. More precisely, we easily have:

\begin{prop} \label{P:doubling} 
If $\bt$ has finite defect $k$, then $\varphi_d(\bt)$ has defect $2k$. More precisely, if $p_1$, \ldots, $p_k$ are the $k$ defective positions in $\bt$, then the defective positions in $\varphi_d(\bt)$ are $2p_i -1$, $2p_i$ for $1 \leq i \leq k$. \qed
\end{prop}

\begin{exa} The periodic infinite word $\bt = (a^2bacaba^2c)^\omega$ has only 2 defects at positions $10$ and $11$, and $\varphi(\bt) = (a^4b^2a^2c^2a^2b^2a^4c^2)^\omega$ 
has $4$ defects at positions $19, 20, 21, 22$. 
\end{exa}

A simple example of a morphism that does not preserve almost richness is $\varphi: a\mapsto ac, b \mapsto b, c \mapsto c$. For instance, consider the (rich) {\em Fibonacci word} $\bbf$, which is generated by the morphism: $a \mapsto ab, b\mapsto a$. We easily see that the image of $\bbf$  by $\varphi$ has only six unique palindromic factors ($\empt, a, b, c, aca, cac$), and hence $\varphi(\bbf)$ has  infinite defect.

\subsection{Class $P$ morphisms}

We now slightly extend the definition of ``class $P$''  morphisms introduced by Hof, Knill, and Simon \cite{aHoKbS95sing} (see also \cite{jAmBjCdD03pali}). 

\begin{defi}[Class $P$ morphisms] \label{D:classP} ~
\begin{itemize} 

\item[i)] A morphism $\varphi$ on $\cA$ is said to be a {\em standard morphism of class $P$} (or a {\em standard $P$-morphism}) if there exists a palindrome $p$ (possibly empty) such that, for all $x\in \cA$, $\varphi(x) = pq_x$ where the $q_x$ are palindromes. If $p$ is non-empty, then some (or all) of the palindromes $q_x$ may be empty or may even take the form $q_x = \pi_x^{-1}$ with $\pi_x$ a proper palindromic suffix of $p$. 

\item[ii)] A morphism $\psi$ on $\cA$ is said to be a {\em morphism of class $P$} (or a {\em $P$-morphism}) if there exists a standard $P$-morphism $\varphi$, with $\varphi(x) = pq_x$ for all $x\in \cA$, such that, for some factorization $p = p'p''$, we have $\psi(x) = p''q_xp'$ for all $x \in \cA$. That is, $\psi = \TT^i(\varphi)$ for some $0 \leq i \leq |p|$.

\end{itemize}
\end{defi}

\begin{remq}
Part $ii)$ of Definition~\ref{D:classP} tells us that any $P$-morphism is a {\em conjugate} of a standard one. Let us also observe that any $P$-morphism as defined in part $ii)$ may also be  a standard $P$-morphism, or a ``dual'' of a standard $P$-morphism (of the form $x \mapsto q_x p$) for other $p$ and $q_x$, because for instance if $|p'| \leq |p''|$, then  $p''q_xp'= (p'' \tilde p'^{-1})( \tilde p' q_x p')$ which has form $r m_x$, where $r$, $m_x$ are palindromes. Indeed, the interest of part $ii)$ is  mainly in view of Definition~\ref{D:special-morphism} hereafter.       
\end{remq}

\begin{note}
The class of $P$-morphisms (resp.~standard $P$-morphisms) is closed under composition, i.e., it is a monoid of morphisms.  
\end{note}

For our purposes, it suffices to consider standard $P$-morphisms in view of the following trivial property.

\begin{prop} \label{P:conjugates}   
Suppose $\varphi$ is a standard $P$-morphism with $\varphi(x) = pq_x$ for all $x\in \cA$ and let $\psi = \TT^i(\varphi)$ for some $i$, $0 \leq i \leq |p|$. Then, for any recurrent infinite word $\bt$, $\psi(\bt)$ and $\varphi(\bt)$ have the same set of factors. \qed
\end{prop}

\begin{exa} The morphism $\tau : a \mapsto baa$, $b\mapsto baba$ is standard $P$ (and its first conjugate $\TT(\tau): a \mapsto aab$, $b\mapsto abab$ is of class $P$). It generates a rich infinite word as does $\TT(\tau)$. This follows easily from the fact that $\tau= \varphi_1 \circ \varphi_2$ with $\varphi_1: a\mapsto a, b\mapsto ba$ and $\varphi_2: a\mapsto ba, b\mapsto bb$, where the latter two morphisms preserve richness: the first one is episturmian and the second one is an {\em insertion morphism} (see Corollary~\ref{C:epi-strict} and Proposition~\ref {P:insert}).
\end{exa}

\begin{defi}  \label{D:special-morphism} We say that a standard  $P$-morphism $\sigma$ is {\em special} if: 1) all $ \sigma(x)=pq_x $ end with different letters, and 2) whenever $\sigma(x)p = pq_x p$, with $x \in \cA$, occurs in some $\sigma(y_1y_2 \cdots y_n)p$, then this occurrence is $\sigma(y_m)p$ for some $m$ with $1\leq m \leq n$. A $P$-morphism is  {\em special} if the corresponding standard  $P$-morphism is special.  
\end{defi}

\begin{remq} 
When $p= \empt$, 2) means that the {\em code} $\sigma(\cA)$ is {\em comma-free} (see \cite{jBdP85theo}). Observe also that the elementary {\em epistandard morphisms} $\{\psi_x \mid x \in \cA\}$ satisfy this definition. Moreover,  as the monoid of epistandard morphisms is generated by all the $\psi_x$ and permutations on $\cA$ (see \cite{xDjJgP01epis, jJgP02epis}), any such morphism is a special  $P$-morphism. For example, $\psi_a \circ \psi_b$ is the special (standard) $P$-morphism with $p = aba$, $q_a = \empt$, $q_b = a^{-1}$.  
\end{remq}

\begin{theo} \label{T: special rich} Suppose $\sigma $ is a special standard $P$-morphism and let $\bt=x_1x_2x_3 \cdots$  be a rich infinite word. Let $h$ be minimal such that  all palindromic factors of $\bt$ of length at most $2$ occur in the prefix $\bt_h$. Then $\sigma (\bt)$ is rich if (and only if) $\sigma (\bt_h)p $ is rich.
\end{theo}
\begin{proof}
 By induction, we suppose $\sigma (\bt_{n-1})p$ is rich for some $n >h$ and show that
$\sigma (\bt_{n})p = \sigma (\bt_{n-1})pq_{x_n}p$ is rich. Let $r$ be the ups of $\bt_n$. Then $R=\sigma (r)p$ is ups of $\sigma (\bt_{n})p$. Indeed, if $R$ has another occurrence in $\sigma (\bt_{n})p$, then by Definition~\ref{D:special-morphism} this occurrence is  $\sigma (x_i\cdots x_j)p$ with $x_i \cdots x_j=r$ and $1\leq i \leq j < n$. This implies that $r$ has another occurrence in $\bt_n$, a contradiction. We have also to show that for any factorization $ef= q_{x_n}p$ with $e,f \neq \empt$, $\sigma (\bt_{n})pf^{-1}$ has a ups. With $r= x_nr'x_n$, $\sigma (\bt_{n})pf^{-1}$ has a palindromic suffix $R'= \tilde f^{-1} Rf^{-1} =\tilde e R'e $. Clearly $r' \neq \empt$, thus if $R'$ has another occurrence in  $\sigma (\bt_{n})pf^{-1}$ then it is 
$\tilde e \sigma (x_i \cdots x_j)pe$. As $e \neq \empt$, we have $x_{i-1}=x_{j+1} =x_n$ and $x_{i-1} \cdots x_{j+1}=r$, a contradiction.  
\end{proof}

\begin{cor} \label{Cor: special rich} Suppose $\sigma$ is a special standard $P$-morphism prolongable on $a$ and let $\bs_k$ be the shortest prefix of $\bs =\sigma^\omega(a)$ that contains all palindromic factors of $\bs$ of length at most $2$. Then $\bs$ is rich if (and only if) $\sigma(\bs_k)p$ is rich.  \qed
\end{cor}

This can be extended to defective words.

 \begin{theo} \label{T: special def} Let $\sigma $ be a special standard $P$-morphism and $\bt$ be an infinite word with finite defect $k$. Let $h$ be minimal such that the prefix $\bt_h$ has defect $k$ and all palindromic factors of $\bt$ of length at most $2$ occur in $\bt_h$. Then $\sigma(\bt)$ is almost rich and its defect is equal to that of $\sigma(\bt_h)p$.
\end{theo}
\begin{proof} Clearly all prefixes $\bt_n$ of $\bt$ with $n>h$ have a ups of length at least $3$. Thus, as in the proof of Theorem~\ref{T: special rich}, we find that all prefixes of $\sigma(\bt)$ longer than $\sigma(\bt_h)p$ have a ups.
\end{proof}

\begin{remq} Naturally one might suspect that if $\sigma$ is a special $P$-morphism prolongable on $a$,  then $\sigma^\omega(a)$ is almost rich. This is not true, as the following proposition shows.
\end{remq}

\begin{prop} The special $P$-morphism $\sigma$: $a \mapsto aba$, $b \mapsto bcb$, $ c \mapsto cac$ generates $\bs = ababcbaba \cdots$ which has infinitely many defects.
\end{prop}
\begin{proof}
Let $p_n= \sigma^n(a)$ and let $w_n$ be the prefix of $\bs$ of length $(3^n+1)/2$, i.e., $w_n= \bs_{(3^n+1)/2}$. Then $w_n$ ends with some letter, $x$ say, which is in the middle of $p_n$. We show by induction that $x$ is the one palindromic suffix of $w_n$. Easily $w_{n+1}x = \sigma(w_n)$, thus $ w_{n+1}$ ends with $y$ such that $xyx=\sigma(x)$. If $ w_{n+1}$ has a palindromic suffix $q$ other than $y$, then easily $|q| > 4$. So it follows by 2) of Definition~\ref{D:special-morphism} that $q= yx\sigma(u)xy$ for some factor $u$ of $\bs$. Hence $\sigma(xux)$ is a palindromic suffix of $ w_{n+1}x$, and therefore  $xux$ is a palindromic suffix of $w_n$, contradicting the induction hypothesis.   
\end{proof}

Indeed we have more generally:

\begin{prop} \label{P:special-generate} Suppose $\sigma$ is a special standard $P$-morphism prolongable on $a$ and let $\bs_h$ be the shortest prefix of $\bs = \sigma^\omega(a)$ that contains all palindromic factors of $\bs$ of length at most $2$. Then $\bs$ has infinite defect if and only if $\sigma(\bs_h)p$ is not rich.
\end{prop}

\begin{proof} ONLY IF: If $\bs$ has infinite defect, then $\sigma(\bs_h)p$ is not rich; otherwise, by Corollary~\ref{Cor: special rich}, $\bs$ would be rich, which is a contradiction. 

IF:  Clearly $\bs$ has at least one defect as $\sigma(\bs_h)p$ is not rich. To show that $\bs$ has infinitely many defects, we suppose by way of contradiction that $\bs$ has finite defect $k \geq 1$. Let $\bs_m$ be the shortest prefix of $\bs$ that has defect $k$. By the minimality of $m$, $\bs_n$ has a ups for all $n \geq m$ and $\bs_m = x_1x_2\cdots x_m$ does not have a ups. But the latter implies that $\sigma(\bs_m)p$ does not have a ups. Indeed, if $\sigma(\bs_m)p$ has a ups, $R$ say, then $R$ begins and ends with $\sigma(x_m)p$. Moreover, as  $\sigma$ is injective, $R = \sigma(x_i \cdots x_m)p$ for some $i \le m$ where $r = x_i\cdots x_m$ is a palindromic suffix of $\bs_m$. But then $r$ must be unioccurrent in $\bs_m$, otherwise $R$ is not unioccurrent in $\sigma(\bs_m)p$, a contradiction. Therefore $\sigma(\bs_m)p$ does not have a ups (i.e., $\bs$ has a defect at position $|\sigma(\bs_m)p| > m$), a contradiction.
\end{proof}

\begin{exa} \label{exemple} Consider the special standard $P$-morphism $\varphi: a \mapsto aab^2aa, b \mapsto bab$.  By Proposition~\ref{P:special-generate},  the infinite words $\varphi^\omega(a)$ and $\varphi^\omega(b)$  have infinitely many defects since their respective prefixes $\varphi(aabb) = aabbaaaabbaababbab$ and $\varphi(babaabb) = babaabbaababaabbaaaabbaababbab$ are not rich (defects at the two penultimate positions in each case). However, if we consider for instance the (rich) Fibonacci word $\bbf$, then $\varphi(\bbf)$ is a rich infinite word. To show this, we need only use Theorem~\ref{T: special rich}: the shortest prefix of $\bbf$ containing all of its palindromic factors of length at most $2$ is $abaa$ and $\varphi(abaa) = aabbaababaabbaaaabbaa$ is rich; whence $\varphi(\bbf)$ is rich. This provides a good example of a non-periodic rich infinite word that is different from a Sturmian word. It was inspired by the family of rich periodic words: $(aab^kaabab)^\omega$ with $k\geq 0$, given in \cite{sBsHmNcR04onth}.
\end{exa}

\begin{remq} From Proposition~\ref{P:special-generate}, we see that special $P$-morphisms generate either rich infinite words or infinite words with infinitely many defective positions. Moreover, as any (primitive) special $P$-morphism generates a uniformly recurrent infinite word with infinitely many palindromic factors, those with infinite defect also have infinitely many oddities (by Proposition~\ref{P:infinite-oddities}).
\end{remq}

\begin{exa} \label{ex:gen-rich} 
Using Corollary~\ref{Cor: special rich}, one can easily verify that the following special standard $P$-morphism generates a rich infinite word: $a \mapsto abb$, $b\mapsto ac$, $c \mapsto a$.
\end{exa}

There is a kind of converse for Theorems \ref{T: special rich} and \ref{T: special def} ({\em cf.} Proposition~\ref{P:defective}).

\begin{theo} \label{T:special-defect} Suppose $\bs = \varphi(\bt)$ where $\varphi$ is a special standard $P$-morphism. Then $D(\bs) \ge D(\bt)$; in particular, if $\bs$ is rich, then $\bt$ is rich.
\end{theo}
\begin{proof}
It suffices to show that if $\bt = x_1x_2x_3 \cdots$ has a defect at position $n$, then $\bs$ has a defect at position $h= |\varphi(\bt_n)p|$. Otherwise, $\bs_h =\varphi(\bt_n)p $ has a ups $R$ beginning and ending with $\varphi(x_n)p$. Thus, as  $\varphi$ is a special $P$-morphism, $R = \varphi(x_i \cdots x_n)p$ for some $i \le n$ where $r = x_i\cdots x_n$ is a palindromic suffix of $\bt_n$. Now, $r$ must be unioccurrent in $\bt_n$, otherwise $R$ is not unioccurrent in $\bs_h$, a contradiction. 
\end{proof}

\begin{rem} \label{R:special-defect}
Notice that property 1) in Definition \ref{D:special-morphism} is too strong here; it suffices that $\varphi$ is injective, i.e., $\varphi(\cA)$ is a code.
\end{rem}

From Theorems~\ref{T: special def} and \ref{T:special-defect}, we immediately see that special $P$-morphisms strictly preserve almost richness. That is:

\begin{theo} Suppose $\bs= \sigma(\bt)$ with $\sigma $ a special $P$-morphism. Then $\bs$ is almost rich if and only if $\bt$ is almost rich. \qed
\end{theo}

Using the following easy lemmas, some of which are well-known, we end this section by proving a theorem which brings us one step closer to a characterization of morphisms preserving richness.

\begin{lem} \label{l:pq} If $p,\ q,\ p',\ q'$ are non-empty palindromes and $pq=p'q'$ is primitive, then $p=p',\ q=q'$. If $pq$ is a primitive palindrome with $p$, $q$ palindromes, then $p$ or $q$ is empty.
\end{lem}
 
\begin{lem} \label{l:pqr} If $pqr$ is a palindrome with $p$, $q$, $r$ palindromes, then $(pq)^h= (rq)^k$, for some $h,k \in \NN, (h,k) \neq (0,0)$.
\end{lem}
 
\begin{lem} \label{l:pal pref} If $Xqp$ is a prefix of $(pq)^\omega$, $pq$ primitive, $p$, $q$ palindromes, then $X = (pq)^ h p$, for some $h \ge 0$.
\end{lem}

\begin{lem} \label{l:expanding} The morphism $\theta: a \mapsto a^n$, $x \mapsto x$ for all letter $x \neq a$ strictly preserves richness. 
\end{lem}

\begin{theo} \label{T:mor pres rich} 
Suppose $\varphi$ is a non-erasing morphism on $\cA$ such that:
\begin{itemize}
\item $\varphi(x) \ne \varphi(y)$ for all letters $x \ne y$;
\item $\varphi(x)$ is a primitive word for any letter $x \in \cA$;
\item  for any three distinct letters $a,b, c$,
\begin{equation}\varphi(a)^\alpha  \varphi(b)^\beta \varphi(c)^\gamma = \empt,\ \alpha , \beta, \gamma \in \ZZ\ \Rightarrow \alpha \beta \gamma=0. \label{e:condi}\end{equation}
\end{itemize}
Then if $\varphi$ preserves richness, it is of class $P$.
\end{theo}
\begin{proof}
Let us denote the images of the letters by $\varphi_i$, $1 \leq  i \leq |\cA|$. We first show that $\varphi(a)=\varphi_1$ and $ \varphi(b)=\varphi_2$ have the form given by the definition of class $P$. As $a^\omega$ is rich, $\varphi_1(a)^\omega$ is rich, so by Theorem~\ref{T:periodic}  $\varphi_1= p_1q_1$ with $p_1$, $q_1$ palindromes. Similarly $\varphi_2= p_2q_2$ with $p_2$, $q_2$ palindromes. Now, as $(a^m b a^m) ^ \omega$ is rich for any $m$, by the same argument as above, we have $\varphi_1^m \varphi_2 \varphi_1^m= PQ$ for some palindromes $P,\ Q$. We shall suppose first that both $\varphi_1 $ and $\varphi_2 $ are not  palindromes. There are three cases according to the place of the separation between $P$ and $Q$.
\begin{itemize}
\item Case $P =\varphi_1^m  X$, $Q= Y \varphi_1^m$, $XY= \varphi_2$, $X,Y \in \cAstar$. If $m$ is large, $p_1q_1X$ is a suffix of $P$ , thus $\tilde X q_1 p_1$ is a prefix of $\varphi_1^m$. By Lemma \ref{l:pal pref}
$\tilde X = p_1(q_1p_1)^\alpha$. Similarly, $Y= (q_1p_1)^\beta q_1$. Thus $\varphi_2=p_1(q_1p_1)^{\alpha + \beta} q_1= \varphi_1^{\alpha + \beta +1} $  which is impossible.
\item Case $PX =\varphi_1^m$, $Q= X \varphi_2 \varphi_1 ^m$, $X \in \cAstar$.
Thus $X \varphi_2 P X =Q$; whence $X$ is a palindrome and  also  $ \varphi_2 P$, i.e., $p_2q_2 P$. By Lemma \ref{l:pqr} and as $p_2q_2$ is primitive, we have $Pq_2=(p_2q_2)^{(\mu+1)}$, and therefore $P= (p_2q_2)^ \mu p_2$. Consider two subcases.
 \begin{itemize}
 \item Case $|P| \ge |p_1q_1|$. As $P$ is a (palindromic) prefix of $\varphi_1^m$ ending with $q_1 p_1$, it has the form $(p_1q_1)^ \lambda p_1$, whence
\begin{equation} (p_1q_1)^ \lambda p_1= (p_2q_2)^ \mu p_2, \ \lambda,\ \mu \ge 0 \label{e:pal} \end{equation}.  
 \item Case $|P| < |p_1q_1|$. In this case, $|X|$ is large and, as it is a palindromic suffix of $\varphi_1^m$,  $X=(q_1p_1)^ \alpha q_1$. Thus, since  $PX =\varphi_1^m$, we get $P= (p_1q_1)^{m-\alpha -1}p_1$, and hence $P=p_1$. So we again get equation \eqref{e:pal} with $\lambda=0$.  
 \end{itemize}

 Now let $p=(p_1q_1)^ \lambda p_1= (p_2q_2)^ \mu p_2$. Then $\varphi_1= p_1q_1=p q_a$ with $q_a= ((p_1q_1)^{\lambda-1}p_1)^{-1}$ and $\varphi_2= p_2q_2=p q_b$ with $q_b= ((p_2q_2)^{\mu-1}p_2)^{-1}$.

\item Case $P =\varphi_1^m \varphi_2 X$, $XQ= \varphi_1^m$, $X,\ Y \in \cAstar$. By symmetry we get 
\begin{equation} (q_1p_1)^ \lambda q_1= (q_2p_2)^ \mu q_2 =p, \ \lambda,\ \mu \ge 0, \label{e:sym} \end{equation}
 and $\varphi_1=q_a p$, $\varphi_2=q_b p$. 
\end{itemize}
 
Now suppose for instance that $\varphi_1$ is a palindrome (but not $\varphi_2$). Then it is easily seen that equation \eqref{e:pal} (resp.~equation~\eqref{e:sym}) also holds with $p_1= \varphi_1$, $q_1= \empt$ (resp.~$p_1= \empt$, $q_1= \varphi_1$). 
Let us also observe that the pair $( \lambda,\ \mu)$ in equation \eqref{e:pal} or \eqref{e:sym} is unique. Indeed if $ \lambda ' > \lambda$ and $\mu ' > \mu$ also work, we get $\varphi_1^{\lambda'-\lambda}= \varphi_2^{\mu'-\mu}$; whence easily  $\varphi_1 =\varphi_2$, a contradiction.
 
Thus the `shape' of class $P$ is satisfied for letters $a,\ b$, $a\ne b$. It remains to pass to $\cA$ in totality. Suppose first, with notations as before, that $\cA$ contains at least two different letters, $a$, $b$ with  $\varphi_1 \neq \varphi_2$ both non-palindromes. Let $c$ be any other letter and $\varphi(c)=\varphi_3=p_3q_3$. Consider the three pairs of letters.
\begin{itemize}
\item First case:
 
Using $(a,b)$: $\varphi_1=p q_a$, $\varphi_2=p q_b$;
 
Using $(a,c)$: $\varphi_1=r s_a$, $\varphi_3=r s_c$;
 
Using $(b,c)$: $\varphi_2=t u_b$, $\varphi_3=t u_c$.
 
Here, $p$, $q$, $r$ are given by equation \eqref{e:pal} and similar ones. Suppose for instance $|p|\ge |r| \ge |t|$. Then we get $p= \varphi_1^\alpha r$, $r= \varphi_3^\beta t$  $p= \varphi_2^\gamma t$ for some $\alpha$, $\beta$, $\gamma$. Thus, $\varphi_1^\alpha \varphi_3^\beta t= \varphi_2^\gamma t$. This gives  $\alpha \beta \gamma =0$ by condition \eqref{e:condi}; whence $p=r$ or $r=t$ or $p=t$. The case $r=t$ for instance gives $\varphi_1=rs_a$, $\varphi_2=r u_b$, $\varphi_3=r s_c$. But, by the observation above, $r=p$.
   
\item Second case: the same for $(a,b)$ and $(a,c)$, but $(b,c)$ gives $\varphi_2= u_b t$, $\varphi_3=u_c t$ and $t$ is given by an equation of form (\ref{e:sym}). We deduce $p = \varphi_1 ^\xi r$,  $pt = \varphi_2 ^\eta $, $rt =  \varphi_3 ^\tau $ for some $\xi$, $\eta$, $ \tau$; whence $ \varphi_2 ^\eta= \varphi_1 ^\xi  \varphi_3 ^\tau $. Clearly, $pt,\ rt \neq \empt$. Thus by (\ref{e:condi}), $\xi=0$, $r=p$.
\end{itemize}
In conclusion, $\varphi$ is a $P$-morphism.

 Now suppose $\cA$ contains exactly one letter, $a$, whose image $\varphi_1$ is not a palindrome and consider any other two letters $b$, $c$.

\begin{itemize}
\item First case:
 
Using $(a,b)$: $\varphi_1=p q_a$, $\varphi_2=p q_b$;
 
Using $(a,c)$: $\varphi_1=r s_a$, $\varphi_3=r s_c$.

Here, $p$, $r$, are given by equation \eqref{e:pal} and a similar one. Suppose for instance $|p|\ge |r|$. We have $p=(p_1q_1) ^ \lambda p_1= \varphi_2 ^\mu$ and $r=(p_1q_1) ^ {\lambda '} p_1= \varphi_2 ^{\mu'}$; whence $\varphi_1^{\lambda ' - \lambda} = \varphi _2^{\mu ' - \mu}$. As $\varphi_1 \neq \varphi_2$ are primitive,  $\lambda' - \lambda =0$, $r=p$. 
\item Second case: the same for $(a,b)$, but $(a,c)$ gives $\varphi_1= s_a r$, $\varphi_3= s_c r$ with $r=(q_1 p_1)^\theta q_1= \varphi_3^\nu$, and hence $ pr=\varphi_1 ^{\lambda + \theta +1}= \varphi_2 ^\mu \varphi_3 ^\nu$. As $\lambda + \theta +1 >0$ this gives $\mu \nu =0$ by \eqref{e:condi}, which is impossible.
\end{itemize}
 
Lastly, if all images of letters are palindromes, then $\varphi$ is trivially of class~$P$.
\end{proof}

\begin{rem} Condition~\eqref{e:condi} of Theorem~\ref{T:mor pres rich}  is satisfied if $ \varphi$ is injective, or if it strictly preserves richness  (using the property that $a^x b^y c ^z a^x$ is not rich). The theorem could be extended to non-primitive $ \varphi(x)$ using Lemma~\ref{l:expanding} but conditions should be formulated accordingly.
\end{rem}

Now let us recall from Theorem~\ref{T:per-alm} that periodic almost rich words are of the form $u ^ \omega$ where $u$ is a product of two palindromes and that only this property is used in proof of Theorem~\ref{T:mor pres rich}. Thus we also have:
      
\begin{theo}Suppose $\varphi$ is a morphism satisfying the conditions of Theorem~$\ref{T:mor pres rich}$. Then if $\varphi$ preserves almost richness, it is of class $P$. \qed
\end{theo}

Furthermore, it is not too difficult to see that `preserves almost richness' could be replaced by `preserves infiniteness of palindromic factors'. This is related to the following long-standing open question posed by Hof, Knill, and Simon in \cite{aHoKbS95sing}: are there (uniformly recurrent) infinite words containing arbitrarily long palindromes that arise from primitive morphisms, none of which belongs to class $P$? The answer is believed to be no. Up to now, it has only been shown to hold in the periodic case (see \cite{jAmBjCdD03pali}) and also in the $2$-letter case (see \cite{bT07mirr}).

\section{Concluding remarks}

To end, we mention a particularly relevant result that gives a good estimate of the palindromic complexity of uniformly recurrent infinite words in terms of the factor complexity. Let us first recall that the {\em palindromic complexity} function $\cP(n)$ (resp.~{\em factor complexity} function $\cC(n)$) of a given infinite word counts the number of different palindromic factors (resp.~number of different factors) of length $n$ for each $n\geq0$. In \cite{pBzMeP07fact},  Bal\'a\v{z}i {\em et al.} proved that for uniformly recurrent infinite words with factors closed under reversal, 
\begin{equation} \label{eq:pali-bound}
  \cP(n) + \cP(n+1) \leq \cC(n+1) - \cC(n) + 2 \quad \mbox{for all $n \in \NN$}.
 \end{equation}

Infinite words for which $\cP(n) + \cP(n+1)$ always reaches the upper bound given in relation~\eqref{eq:pali-bound} can be viewed as words containing the maximum number of palindromes. Naturally one would conjecture that all such words are rich.  Indeed, this assertion is true -- it was recently proved by the first and fourth authors together with M.~Bucci and A.~De Luca in \cite{mBaDaGlZ08acon}.  Interestingly, its proof relies upon another new characterization of rich words, which is useful for establishing the key part of the proof, namely that the so-called {\em super reduced Rauzy graph} is a tree. 


\begin{thebibliography}{99}  


\bibitem{jAmBjCdD03pali} J.-P.~Allouche, M.~Baake, J.~Cassaigne, D.~Damanik, Palindrome complexity, {\it Theoret. Comput. Sci.} {292} (2003) 9--31.

\bibitem{pAzMePcF06pali} P.~Ambro{\v{z}}, C.~Frougny, Z.~Mas{\'a}kov{\'a}, E.~Pelantov{\'a},  Palindromic complexity of infinite words associated with 
  simple {P}arry numbers, {\it Ann. Inst. Fourier (Grenoble)} 56 (2006) 2131--2160.
  
 \bibitem{vAlZiZ03pali} V.~Anne, L.Q.~Zamboni, I.~Zorca, Palindromes and
  pseudo-palindromes in episturmian and pseudo-palindromic infinite words, in: {\it Proceedings of the Fifth International Conference on Words} (Montr\'eal, Canada), September 13--17, 2005. {\it Publications du LaCIM} 36 (2005) 91--100.
 

\bibitem{pBzMeP07fact} P.~Bal\'a\v{z}i, Z.~Mas\'akov\'a, E.~Pelantov\'a, Factor versus
  palindromic complexity of uniformly recurrent infinite words, {\em Theoret. 
  Comput. Sci.} {380} (2007) 266--275. 
  
  
  \bibitem{jB07stur}
J.~Berstel, Sturmian and episturmian words (A survey of some recent results), in: {\em Proceedings of CAI 2007}, Lecture Notes in Computer Science, vol.~4728, 2007,  pp.~23--47.

\bibitem{jBdP85theo} J.~Berstel, D.~Perrin, \emph{Theory of codes}, vol.~117 of {\em Pure and Applied
  Mathematics}, Academic Press Inc., Orlando, FL, 1985. 

\bibitem{sBsHmNcR04onth} S.~Brlek, S.~Hamel, M.~Nivat, C.~Reutenauer,  On the palindromic complexity of infinite words, {\it Internat. J. Found. Comput. Sci.} { 15} (2004) 293--306.

\bibitem{mBaDaGlZ08acon} M.~Bucci, A.~De~Luca, A.~Glen, L.Q.~Zamboni, A connection between palindromic and factor complexity using return words, {\it Adv. in Appl. Math.}, to appear, arXiv:0802.1332.

\bibitem{jC97comp} J.~Cassaigne, Complexit\'e et facteurs sp\'eciaux,  {\it Bull. Belg. Math. Soc. Simon Stevin} { 4} (1997) 67--88.

\bibitem{eCgH73sequ} E.M.~Coven, G.A.~Hedlund, Sequences with minimal block growth, \emph{Math. Systems Theory} { 7} (1973) 138--153.

\bibitem{dDlZ03comb}
D.~Damanik, L.Q.~Zamboni, Combinatorial properties of {A}rnoux-{R}auzy subshifts and applications to {S}chr\"odinger operators, {\it Rev. Math. Phys.} 15 (2003) 745--763.


\bibitem{aD97stur} A.~de~Luca, Sturmian words: structure,
combinatorics and their arithmetics, \emph{Theoret. Comput. Sci.}
{183} (1997) 45--82. 

\bibitem{xDjJgP01epis} X.~Droubay, J.~Justin, G.~Pirillo,
Episturmian words and some constructions of de Luca and Rauzy,
\emph{Theoret. Comput. Sci.} {255} (2001) 539--553.  


\bibitem{xDgP99pali}
X.~Droubay, G.~Pirillo, Palindromes and {S}turmian words, {\it Theoret. 
  Comput. Sci.} 223 (1999) 73--85. 

  \bibitem{sF06pali}
S.~Fischler, Palindromic prefixes and episturmian words,
{\it J. Combin.  Theory Ser.~A} 113 (2006) 1281--1304.

\bibitem{sF06pali2}
S.~Fischler, Palindromic prefixes and diophantine approximation, {\it Monatsh. Math.} 151 (2007) 11--37.

\bibitem{aF73comp} A.S.~Fraenkel, Complementing and exactly covering sequences, {\it J. Combin.  Theory Ser. A} { 14} (1973) 8--20.    

\bibitem{aGjJ07epis} A.~Glen, J.~Justin, Episturmian words: a survey, Preprint, 2007, arXiv:0801.1655.

\bibitem{aHoKbS95sing} A.~Hof, O.~Knill, B.~Simon, Singular continuous spectrum for palindromic 
Schr\"{o}dinger operators, {\it Commun. Math. Phys.} {174} (1995) 149--159.

\bibitem{cHlZ99desc} C.~Holton, L.Q.~Zamboni, Descendants of primitive substitutions, {\em Theory Comput. Syst.} { 32} (1999) 133--157. 


\bibitem{jJgP02epis} J.~Justin, G.~Pirillo, Episturmian words and episturmian morphisms, \emph{Theoret. Comput. Sci.} {276} (2002) 281--313. 

\bibitem{jJlV00retu} J.~Justin, L.~Vuillon, Return words in {S}turmian and episturmian words, {\em Theoret. Inform. Appl.} { 34} (2000) 343--356.


\bibitem{mL83comb} M.~Lothaire, {\em Combinatorics On Words}, vol.~17 of {\em Encyclopedia of Mathematics and its Applications}, Addison-Wesley, Reading, Massachusetts, 1983.

\bibitem{mL02alge} M.~Lothaire, \emph{Algebraic Combinatorics On
Words}, vol.~90 of {\em Encyclopedia of Mathematics and its Applications}, Cambridge University Press, U.K., 2002.

\bibitem{gPlV07acha} G.~Paquin, L.~Vuillon, A characterization of balanced episturmian sequences, {\it Electron. J. Combin.} 14 (2007) \#R33, pp.~12. 

\bibitem{nP02subs} N.~Pytheas Fogg, {\em Substitutions In Dynamics, Arithmetics And Combinatorics}, vol.~1794 of {\em Lecture Notes in Mathematics}, Springer-Verlag, Berlin, 2002.

\bibitem{mQ87subs}
M.~Queff\'{e}lec, \emph{Substitution dynamical systems -- spectral analysis}, vol.~1924 of {\it Lecture Notes in Mathematics}, Springer-Verlag, New York, 1987.

\bibitem{bT07mirr}
B.~Tan, Mirror substitutions and palindromic sequences, {\it Theoret. Comput. Sci.} 389 (2007) 118--124.

\bibitem{Tij1} R.~Tijdeman,  Exact covers of balanced sequences and Fraenkel's conjecture, in: {\it Algebraic Number Theory and Diophantine Analysis} (Graz, $1998$), de Gruyter, Berlin, 200, pp.~467--483.


\bibitem{lV03bala} L. Vuillon, Balanced words, {\it Bull. Belg. Math. Soc. Simon Stevin} 10 (2003) 787--805.

\end{thebibliography}
\end{document}